\documentclass[11pt,letterpaper]{amsart}
\usepackage{amsthm,amsfonts,amsmath,amssymb,latexsym}
\usepackage{epsfig,graphics,color}
\usepackage[all]{xy}
\usepackage[breaklinks=true]{hyperref}
\usepackage{mathrsfs}
\usepackage{stmaryrd}
\usepackage{verbatim}
\usepackage{bm}
\usepackage{mathabx}
\usepackage{enumitem}
\usepackage{pgf,amsmath,tikz,type1cm,fix-cm}
\usetikzlibrary{positioning}
\usepackage{xcolor}
\newcommand{\new}[1]{{\color{red}#1}}
\newenvironment{redblock}{\begingroup\color{red}}{\endgroup}

\newenvironment{blueblock}{\begingroup\color{blue}}{\endgroup}

\setlist[itemize,1]{leftmargin=\dimexpr 26pt-.1in}

\newtheorem{PARA}{}[section]
\newtheorem{theorem}[PARA]{Theorem}
\newtheorem{corollary}[PARA]{Corollary}
\newtheorem{lemma}[PARA]{Lemma}
\newtheorem{proposition}[PARA]{Proposition}
\newtheorem{definition}[PARA]{Definition}
\newtheorem{definition-proposition}[PARA]{Definition-Proposition}
\newtheorem{conjecture}[PARA]{Conjecture}
\newtheorem{question}[PARA]{Question}
\theoremstyle{definition}
\newtheorem{remark}[PARA]{Remark}
\newtheorem{example}[PARA]{Example}

\newcommand{\para}{\begin{PARA}\rm}
\newcommand{\arap}{\end{PARA}\rm}
\newcommand{\dfn}{\begin{definition}\rm}
\newcommand{\nfd}{\end{definition}\rm}
\newcommand{\rmk}{\begin{remark}\rm}
\newcommand{\kmr}{\end{remark}\rm}
\newcommand{\xmpl}{\begin{example}\rm}
\newcommand{\lpmx}{\end{example}\rm}

\newcommand{\N}{{\mathbb{N}}}

\newcommand{\R}{{\mathbb{R}}}

\newcommand{\Z}{{\mathbb{Z}}}
\newcommand{\Int}{\mathrm{Int}\,}

\newcommand{\eps}{{\varepsilon}}
\newcommand{\om}{{\omega}}

\def\NABLA#1{{\mathop{\nabla\kern-.5ex\lower1ex\hbox{$#1$}}}}
\def\Nabla#1{\nabla\kern-.5ex{}_{#1}}
\def\Tabla#1{\Tilde\nabla\kern-.5ex{}_{#1}}
\renewcommand{\Tilde}{\widetilde}

\newcommand{\p}{{\partial}}

\newcommand{\ol}{\overline}

\newcommand{\wt}{\widetilde}

\DeclareMathOperator{\cov}{cov}
\DeclareMathOperator{\Edge}{Edge}
\newcommand{\ceil}[1]{\left\lceil #1 \right\rceil}

\DeclareMathOperator{\vol}{vol}
\DeclareMathOperator{\Symp}{Symp}
\begin{document}

\title{Covering numbers of metric spaces}
\author{Kai Cieliebak}
\email{kai.cieliebak@math.uni-augsburg.de}
\author{Milan Cieliebak}
\email{milan.cieliebak@gmail.com}
\date{\today}


\begin{abstract}
We propose and provide evidence for the following conjecture: the
covering number of a product of metric spaces with the sum metric
is greater than or equal to the product of their covering numbers.
This generalizes Vizing's conjecture in graph theory. 
We prove an infinitesimal version of the conjecture for Finsler
manifolds of dimension at least six. 
\end{abstract}

\maketitle

\section{Introduction}\label{sec:introduction}

For a metric space $(X,d)$ and a radius $r>0$ define the {\em covering number}
$$
  \cov(X,d,r) := \inf\{n\in\N_0\mid \text{ $X$ can be covered by $n$
    closed $r$-balls}\}.
$$
If there is no such covering we set $\cov(X,d,r):=\infty$. 
In particular, $\cov(\varnothing)=0$. We use the conventions
$0\cdot\infty=0$ and $n\cdot\infty=\infty$ for $n\geq 1$.
If $d$ and/or $r$ are clear from the context they will be omitted.

On the product $X\times Y$ of two metric spaces $(X,d_X)$, $(Y,d_Y)$
consider the metrics
\begin{equation}\label{eq:dp}
  d_p\bigl((x_1,y_1),(x_2,y_2)\bigr) := \begin{cases}
    \bigl(d_X(x_1,x_2)^p+d_Y(y_1,y_2)^p\bigr)^{1/p} & 1\leq p<\infty,\cr
    \max\{d_X(x_1,x_2),d_Y(y_1,y_2)\} & p=\infty.
  \end{cases}
\end{equation}
It is the goal of this note to investigate to what extent the
following {\em covering inequality} holds: 
\begin{equation}\label{eq:cov}
  \cov(X\times Y,d_p,r) \geq \cov(X,d_X,r)\,\cov(Y,d_Y,r).
\end{equation}
Note that by rescaling it suffices to consider the inequality for
$r=1$, but it will sometimes be convenient to allow other radii.

A simple counterexample (Example~\ref{ex:Z4} below) shows
that~\eqref{eq:cov} fails in general for $p>1$. This leaves us with
the following 

\begin{conjecture}\label{conj:cov}  
For all metric spaces $(X,d_X)$, $(Y,d_Y)$ and each radius $r>0$ we have
$$
  \cov(X\times Y,d_1,r) \geq \cov(X,d_X,r)\,\cov(Y,d_Y,r).
$$
\end{conjecture}

An interesting special case of this is Vizing's conjecture in graph
theory, see e.g.~\cite{Gaar19}.
Given a graph $G$ (which will always be assumed finite and connected)
with vertex set $X$, a {\em dominating set} is a subset $D\subset X$
such that each vertex in $X\setminus D$ is adjacent to some point in $D$. 
The {\em domination number} $\gamma(G)$ of $G$ is the smallest
cardinality of a dominating set. If we equip $X$ with the graph metric
on $G$ with edge length $1$, then with $r=1$ we have 
$$
  \gamma(G) = \cov(X).
$$
Consider now two graphs $G,H$ with vertex sets $X,Y$ equipped with
the graph metrics of edge length $1$ and let $r=1$. Let $G\square H$
be the graph with vertex set $X\times Y$ and edges $[(x,y),(x,y')]$ for
$x\in X$ and $[y,y']\in\Edge(H)$, and $[(x,y),(x',y)]$ for
$[x,x']\in\Edge(G)$ and $y\in Y$. Then
$$
  \gamma(G\square H) = \cov(X\times Y,d_1),
$$
so the following is a special case of Conjecture~\ref{conj:cov}:
  
\begin{conjecture}[Vizing~\cite{Vizing63,Vizing68}]\label{conj:Vizing}  
For all graphs $G,H$ we have
$$
  \gamma(G\square H) \geq \gamma(G)\,\gamma(H).
$$
\end{conjecture}

While we did not succeed in proving Conjecture~\ref{conj:cov}, we will
provide evidence by proving some special and limiting cases.

Our main result is the following infinitesimal version of Conjecture~\ref{conj:cov}.
A (continuous) {\em Finsler metric} on a manifold $M$ is a continuous
fibrewise convex and positively homogeneous function $F:TM\to\R_{\geq
  0}$ such that $F^{-1}(0)$ is the zero section. 
It gives rise to a length of curves, and thus to a (possibly
nonsymmetric) metric $d_M$ on $M$. 
See~\S\ref{sec:Finsler} for details.

\begin{theorem}\label{thm:Finsler-intro}
For each pair of compact Finsler manifolds $M,N$ of dimensions $m\geq n>5$
there exist constants $r_{M,N}>0$ and $p_{m,n}>1$ such that for all
$0<r\leq r_{M,N}$ and $1\leq p\leq p_{m,n}$ we have
$$
  \cov(M\times N,d_p,r) > \cov(M,d_M,r)\,\cov(N,d_N,r).
$$
\end{theorem}

The proof of Theorem~\ref{thm:Finsler-intro} reduces by localization
and approximation arguments to an inequality for covering densities of
convex bodies.
The dimension constraint is due to the lack of good upper bounds on
the covering density in low dimensions. More precisely, the theorem
holds for all pairs $(n,m)$ with $1\leq n\leq m$ 
except possibly the following ones: $(1,m)$ with $m\geq 3$,
$(2,m)$ with $3\leq m\leq 12$, $(3,m)$ with
$3\leq m\leq 15$, $(4,m)$ with $4\leq m\leq 8$, $(5,5)$, and $(5,6)$.

If $M,N$ are Riemannian manifolds, then their product $M\times N$ is
again a Riemannian manifold whose metric is obtained from $d_M,d_N$ by
formula~\eqref{eq:dp} with $p=2$. In this case, by proving a lower
bound on the covering density, we can improve
Theorem~\ref{thm:Finsler-intro} to

\begin{theorem}\label{thm:Riemann-intro}
For each pair of compact Riemannian manifolds $M,N$ of dimensions $m\geq n>2$
there exists a constant $r_{M,N}>0$ such that for all $0<r\leq r_{M,N}$ we have
$$
  \cov(M\times N,d_2,r) > \cov(M,d_M,r)\,\cov(N,d_N,r).
$$
\end{theorem}

More precisely, Theorem~\ref{thm:Riemann-intro} holds for all $m\geq
n\geq 1$ except possibly for $n=2$ with $25\leq m\leq 75$, and $n=1$ with
$16\leq m\leq 4428$. 

{\bf Further results and structure of the paper. }
We begin in~\S\ref{sec:basic} with some basic observations such as the
following one: 
{\em If Conjecture~\ref{conj:cov} holds for all finite metric spaces
  $X,Y$ then it also holds for all compact metric spaces $X,Y$.} 

In~\S\ref{sec:extremal} we introduce the notion of an {\em extremal pair} $(B,b)$,
consisting of an $r$-ball $B\subset X$ and a point $b\in B$ such that
every $r$-ball containing $b$ is contained in $B$.
We show that each finite subset of a metric tree has an extremal pair
and use an inductive procedure to prove: 
{\em Conjecture~\ref{conj:cov} holds if $X$ is a subset of a metric tree
and $Y$ is arbitrary.}
Applied to graphs, this recovers Vizing's conjecture in the case that
one graph is a tree.
The proof uses properties of the Gromov product on $0$-hyperbolic metric spaces.

In~\S\ref{sec:density} we show that an infinitesimal version of
inequality~\eqref{eq:cov} on affine space is equivalent to an
inequality for covering densities. We recall upper bounds on the
covering density from the literature, and we prove a new lower bound
on the covering density of ellipsoids. 

In~\S\ref{sec:Finsler} we use the bounds from~\S\ref{sec:density} to
prove Theorems~\ref{thm:Finsler-intro} and~\ref{thm:Riemann-intro}.
\medskip

{\bf Acknowledgements. }
Many thanks to Elisabeth Gaar for drawing our attention to Vizing's
fascinating conjecture. We also thank
Filip Bro\'ci\'c,
Javier G\'omez-Serrano,
Suny Kim,
Alexandru Oancea and
Milan Zerbin
for interesting discussions. 

{\bf AI disclosure. }
The lower bound on the volume $W_n$ in~\S\ref{ss:Wn-bound} was found
with the help of GPT-5.2-pro, and we used Claude Opus 5.5 for
final corrections. Other than this, no generative AI was used in the
preparation of this note.

\section{Basic observations}\label{sec:basic}

In the following we fix a radius $r>0$ and drop it from the notation.
We will also drop the metrics $d_X,d_Y$ from the notation when they
are clear from the context. 

{\bf Relative covering numbers. }
For a metric space $(X,d)$ and $x\in X$ we denote by
$B_X(x):=\{x'\in X\mid d(x,x')\leq r\}$ the closed $r$-ball in $X$
with center $x$. For subsets $A,Z\subset X$ let
$$
  \cov_Z(A) := \inf\Bigl\{n\in\N_0\;\Big|\; A\subset \bigcup_{i=1}^n B_X(z_i)
  \text{ for some }z_1,\dots,z_n\in Z\Bigr\}
$$
be the minimal number of $r$-balls {\em with centers in $Z$} needed to
cover $A$, with $\cov_Z(A):=\infty$ if there is no such covering.
The {\em relative covering number} of $A\subset X$ is
$$
  \cov(A;X):=\cov_X(A).
$$
By contrast, for a subset $Y\subset X$,
$\cov(Y)$ always denotes the covering number of the metric space $Y$
with the induced metric, i.e., with centers in $Y$; we call it the
{\em intrinsic} covering number. Thus $\cov(X)=\cov(X;X)$, and
$\cov(Y;X)\leq\cov(Y)$ for $Y\subset X$.

\begin{lemma}\label{lem:basic}
(a) For all metric spaces $X,Y$ we have
$$
  \cov(X\times Y,d_\infty) \leq \cov(X)\,\cov(Y).
$$
(b) For a subset $Y\subset (X,d_X)$ and its closure $\ol Y$ we have
$$
  \cov(\ol Y)\leq\cov(Y)\qquad\text{and}\qquad \cov(\ol Y;X)=\cov(Y;X).
$$
(c) For any closed $r$-ball $B\subset X$ we have
$$
  \cov(X\setminus B)\geq\cov(X)-1. 
$$
(d) For $p\leq q$ and all metric spaces $X,Y$ we have
$$
  \cov(X\times Y,d_p) \geq \cov(X\times Y,d_q).
$$
\end{lemma}

\begin{proof}
Part (a) holds because $r$-balls in the metric $d_\infty$ are the
products of $r$-balls, and for coverings of $X$ and $Y$ their products
cover $X\times Y$.
Part (b) holds because a finite union of closed balls is closed. 
Indeed, if the balls $B_X(y_i)$ with $y_i\in Y$, $i=1,\dots,n$, cover
$Y$, then their union is closed and hence also covers $\ol Y$, so the
balls $B_{\ol Y}(y_i)=B_X(y_i)\cap\ol Y$ with centers $y_i\in\ol Y$
cover $\ol Y$. The same argument with centers $x_i\in X$ gives
$\cov(\ol Y;X)\leq\cov(Y;X)$, and the reverse inequality is clear
from $Y\subset\ol Y$.
For (c) note that if the closed $r$-balls $B_1,\dots,B_n$ cover
$X\setminus B$, then the balls $B_1,\dots,B_n,B$ cover $X$. 
Part (d) follows from $d_p\geq d_q$. 
\end{proof}

\begin{example}\label{rem:basic-b}
For $X=\R$, $r=1$ and $Y=(0,1)\cup(1,2)$ we have $\ol Y=[0,2]=B_{\ol Y}(1)$,
so $\cov(\ol Y)=1$, whereas $\cov(Y)=2$ because the only possible center
of a single ball covering $Y$ would be $1\notin Y$.
\end{example}

\begin{lemma}\label{lem:min-covering}
Let $B_1,\dots,B_n$ be a minimal covering of a nonempty space $X$. Then there exists a
point of $X$ that lies in only one of the balls $B_i$. 
\end{lemma}

\begin{proof}
If each point lies in at least two of the balls, then each point of
$B_n$ lies in $B_1\cup\cdots\cup B_{n-1}$, so $B_1,\dots,B_{n-1}$
cover $X$. 
\end{proof}

\begin{corollary}\label{cor:Milan}
If $\cov(Y)>1$ and $1\leq\cov(X)<\infty$, then $\cov(X\times Y,d_\infty)>\cov(X)$. 
\end{corollary}

\begin{proof}
Since $Y\neq\varnothing$, the $X$-factors of a covering of $X\times Y$ by
balls cover $X$, so at least $\cov(X)$ balls are needed.
Suppose we can cover $X\times Y$ by balls $B_i=B_i^X\times B_i^Y$,
$i=1,\dots,n=\cov(X)$. Then $B_1^X,\dots,B_n^X$ is a minimal
covering of $X$, so by Lemma~\ref{lem:min-covering} some $x\in X$
lies in only one ball $B_i^X$. But then $\{x\}\times Y$ must be
contained in $B_i=B_i^X\times B_i^Y$, hence $Y\subset B_i^Y$, and
therefore $\cov(Y)=1$. 
\end{proof}

{\bf Reduction to finite metric spaces. }

\begin{lemma}\label{lem:finite-subset}
Let $(X,d)$ be a compact metric space, $S\subset X$ a finite subset,
and $r>0$. Then there exists a finite subset $T\subset X$ such that
$S\subset T$ and $\cov(T,d,r)=\cov(X,d,r)$. 
\end{lemma}

\begin{proof}
We will drop $d$ from the notation. 
Take a covering of $X$ by balls $B(y_i,r)$, $i=1,\dots,\cov(X,r)$,
and add the $y_i$ to $S$, still denoting the new set $S$. Then
$\cov(T,r)\leq\cov(X,r)$ for each subset $T\subset X$ containing $S$.

There exists a countable dense subset $C\subset X$. (For example, pick
for each $n\in\N$ a finite covering of $X$ by balls of radius $1/n$
and take for $C$ the set of the centers of all these balls.) 
We may assume $S\subset C$. Write $C=\{c_1,c_2,\dots\}$ such that
$S=\{c_1,\dots,c_s\}$ and denote $C_n:=\{c_1,\dots,c_n\}$. We will
show that there exists $n\geq s$ with $\cov(C_n,r)\geq\cov(X,r)$. By
the preceding discussion, $T=C_n$ is then the desired subset. 

Arguing by contradiction, suppose that $\cov(C_n,r)\leq\cov(X,r)-1=:k$
for all $n\geq s$. Thus for each $n\geq s$ there exist indices
$i(n,1),\dots,i(n,k)$ such that the balls $B(c_{i(n,j)},r)$ with
$j=1,\dots,k$ cover $C_n$. 
By compactness of $X$, there exists an increasing sequence
$n_\ell\to\infty$ such that
$c_{i(n_\ell,j)}\to x_j\in X$ as $\ell\to\infty$ for each $j=1,\dots,k$.
We claim that the balls $B(x_j,r)$ with $j=1,\dots,k$ cover $X$. 

To see this, we again argue by contradiction and suppose there exists
$x\in X$ with $d(x,x_j) \geq r+\eps>r$ for all $j=1,\dots,k$. 
Pick $\ell$ so large that the following hold: 
\begin{itemize}
\item $d(c_{i(n_\ell,j)},x_j)<\eps/2$ for all $j=1,\dots,k$ (by convergence);
\item there exists $c_m\in C_{n_\ell}$ with $d(x,c_m)<\eps/2$ (by
  density of $C$).
\end{itemize}
Since the balls $B(c_{i(n_\ell,j)},r)$ with $j=1,\dots,k$ cover $C_{n_\ell}$,
there exists $j\in\{1,\dots,k\}$ with $d(c_m,c_{i(n_\ell,j)})\leq
r$. For this $j$ we get
$$
  d(x,x_j)
  \leq d(x,c_m) + d(c_m,c_{i(n_\ell,j)}) + d(c_{i(n_\ell,j)},x_j) 
  < \eps/2 + r + \eps/2 = r + \eps. 
$$
This contradiction proves the claim. It follows that $\cov(X,r)\leq k
< \cov(X)$, which is again a contradiction and thus proves the lemma.
\end{proof}
  
\begin{proposition}\label{prop:finite-to-compact}
Fix $r>0$ and $p\in[1,\infty]$. If~\eqref{eq:cov} holds for all finite
metric spaces $X,Y$ then it also holds for all compact metric spaces
$X,Y$. 
\end{proposition}

\begin{proof}
Assume that~\eqref{eq:cov} holds for all finite metric spaces and
consider compact metric space $X,Y$. Set $k:=\cov(X\times Y,d_p,r)$
and pick a covering of $X\times Y$ by $d_p$-balls of radius $r$ with
centers $(x_1,y_1),\dots,(x_k,y_k)$. By Lemma~\ref{lem:finite-subset},
there exist finite subsets $T_X\subset X$ such that
$\{x_1,\dots,x_k\}\subset T_X$ and $\cov(T_X,d_X,r)=\cov(X,d_X,r)$, 
and $T_Y\subset Y$ such that
$\{y_1,\dots,y_k\}\subset T_Y$ and $\cov(T_Y,d_Y,r)=\cov(Y,d_Y,r)$.
Then $T_X\times T_Y$ contains $(x_1,y_1),\dots,(x_k,y_k)$ and is
therefore covered by the $d_p$-balls of radius $r$ with these centers, 
so using~\eqref{eq:cov} for the finite metric spaces $T_X,T_Y$ we get
\begin{align*}
  \cov(X\times Y,d_p,r) &= k \geq \cov(T_X\times T_Y,d_p,r) \cr
  &\geq \cov(T_X,d_X,r)\,\cov(T_Y,d_Y,r) \cr
  &= \cov(X,d_X,r)\,\cov(Y,d_Y,r).
\end{align*}
\end{proof}

\section{Extremal pairs}\label{sec:extremal}

Throughout this section, we fix $r>0$ and drop it from the notation.
Let us call a pair $(B,b)$ consisting of an $r$-ball $B\subset X$ and
a point $b\in B$ {\em extremal} if every $r$-ball containing $b$ is
contained in $B$. 

We will also need the following relative version, where $B_X(x)$ and
$\cov_Z(A)$ are as in~\S\ref{sec:basic}. 
Let $A,Z\subset X$. A {\em $Z$-relative extremal pair for $A$} is a
pair $(B,b)$ consisting of a ball $B=B_X(c)$ with center $c\in Z$ and
a point $b\in A\cap B$ such that 
$$
  B_X(z)\cap A\subset B\qquad\text{for every $z\in Z$ with $b\in B_X(z)$.}
$$
For $Z=X$ we call $(B,b)$ a {\em relative extremal pair for $A$}. For
$A=Z=X$ this is just an extremal pair in $X$. 

\begin{lemma}\label{lem:extremal}
If $(B,b)$ is a $Z$-relative extremal pair for $A\subset X$, then
$$
  \cov_Z(A\setminus B) = \cov_Z(A)-1
$$
(with $\infty-1:=\infty$). In particular, if $(B,b)$ is an extremal pair
in $X$, then
$$
  \cov(X\setminus B;X)=\cov(X)-1. 
$$
\end{lemma}

\begin{proof}
Let $c\in Z$ be the center of $B$. If balls $B_X(z_1),\dots,B_X(z_m)$
with $z_i\in Z$ cover $A\setminus B$, then together with $B=B_X(c)$
they cover $A$. Hence $\cov_Z(A)\leq\cov_Z(A\setminus B)+1$, which
proves the lemma if $\cov_Z(A)=\infty$. 
For the converse inequality, consider a covering $B_1,\dots,B_n$ of $A$ by
balls with centers in $Z$, where $n=\cov_Z(A)<\infty$. Since $b\in A$,
one of the balls, say $B_n$, contains $b$. The definition of a
$Z$-relative extremal pair implies $B_n\cap A\subset B$. Thus
$B_n\cap(A\setminus B)=\varnothing$, so $B_1,\dots,B_{n-1}$ must cover
$A\setminus B$ and we have $\cov_Z(A\setminus B) \leq n-1$. (For $n=1$
this means $A\setminus B=\varnothing$.) 
\end{proof}

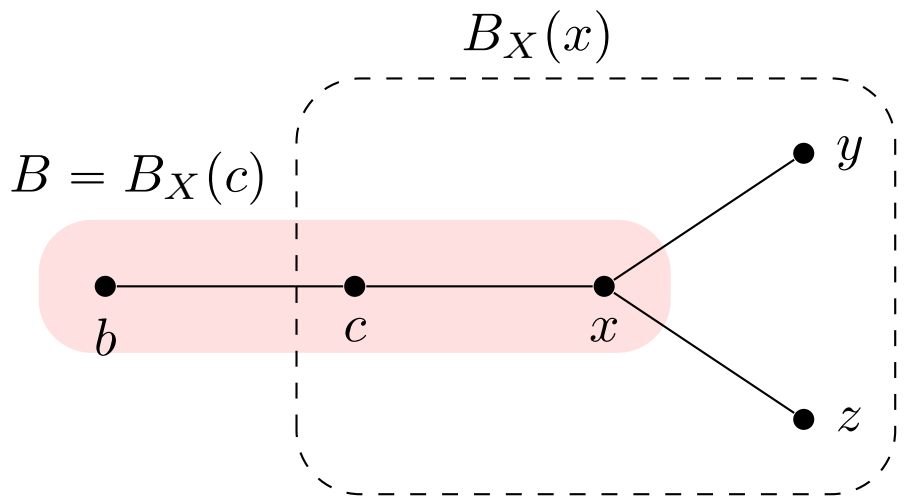
\begin{figure}[ht]
\centering
\begin{tikzpicture}[scale=1.1,every node/.style={font=\small},
    vtx/.style={circle,fill,inner sep=1.4pt}]
  \fill[red!12,rounded corners=10pt] (-0.4,-0.4) rectangle (3.4,0.4);
  \draw[dashed,rounded corners=12pt] (1.15,-1.25) rectangle (4.75,1.25);
  \node[vtx,label=below:$b$] (b) at (0,0) {};
  \node[vtx,label=below:$c$] (c) at (1.5,0) {};
  \node[vtx,label=below:$x$] (x) at (3,0) {};
  \node[vtx,label=right:$y$] (y) at (4.2,0.8) {};
  \node[vtx,label=right:$z$] (z) at (4.2,-0.8) {};
  \draw (b)--(c)--(x)--(y) (x)--(z);
  \node at (0.2,0.65) {$B=B_X(c)$};
  \node at (2.6,1.5) {$B_X(x)$};
\end{tikzpicture}
\caption{The tree in Example~\ref{ex:extremal-intrinsic}.}
\label{fig:extremal-intrinsic}
\end{figure}

\begin{example}\label{ex:extremal-intrinsic}
Lemma~\ref{lem:extremal} fails for intrinsic covering numbers. For
example, let $X=\{b,c,x,y,z\}$ be the set of vertices of the tree
shown in Figure~\ref{fig:extremal-intrinsic} with the graph metric of
edge length $1$, and let
$r=1$. The balls containing $b$ are $B_X(b)=\{b,c\}$ and
$B_X(c)=\{b,c,x\}$, so $(B,b)$ with $B:=B_X(c)$ is an extremal pair. 
We have $\cov(X)=2$ because $B_X(c)\cup B_X(x)=X$ and no single ball
covers $X$. Moreover, $X\setminus B=\{y,z\}\subset B_X(x)$, so
$\cov(X\setminus B;X)=1$ as predicted by Lemma~\ref{lem:extremal}. 
However, $d(y,z)=2>r$, so the intrinsic covering number is
$\cov(X\setminus B)=2=\cov(X)$. The ball $B_X(x)$ covering $X\setminus B$ 
has its center $x$ in $B$. 
\end{example}

Let $Z\subset X$. We call a subset $A\subset X$ {\em $Z$-eliminable} if
there exists a finite chain
$$
  A=A_0\supsetneq A_1\supsetneq\cdots\supsetneq A_k=\varnothing,\qquad
  A_{i+1}=A_i\setminus B_i,
$$
where $(B_i,b_i)$ is a $Z$-relative extremal pair for $A_i$ for
$i=0,\dots,k-1$. We call this a {\em $Z$-eliminating chain of length
$k$}. Since $A_i\setminus A_{i+1}\subset B_i$, the balls
$B_0,\dots,B_{k-1}$ cover $A$, so $\cov_Z(A)\leq k$. By
Lemma~\ref{lem:extremal}, in fact $k=\cov_Z(A)$.

For metric spaces $X,Y$ and $p\in[1,\infty]$ we consider $d_p$-balls
of radius $r$ in $X\times Y$. Since $d_p\geq d_\infty$, the $d_p$-ball
with center $(x,y)\in X\times Y$ is contained in $B_X(x)\times B_Y(y)$,
with equality for $p=\infty$. For $A\subset X$ we denote by
$\cov(A\times Y;X\times Y,d_p)$ the minimal number of $d_p$-balls
with centers in $X\times Y$ that cover $A\times Y$.

\begin{proposition}\label{prop:extremal}
Fix any $p\in[1,\infty]$ and $r>0$. 
Let $Z\subset X$, let $A\subset X$ have a $Z$-eliminating chain of
length $k$, and let $Y$ be any metric space. If $N$ $d_p$-balls with
centers in $Z\times Y$ cover $A\times Y$, then
$$
  N\geq k\,\cov(Y) = \cov_Z(A)\,\cov(Y).
$$
(In particular, no such finite covering exists if $k\geq 1$ and
$\cov(Y)=\infty$.) Hence, if $A$ is $X$-eliminable, then 
$$
  \cov(A\times Y;X\times Y,d_p) \geq \cov(A;X)\,\cov(Y),
$$  
and if $X$ is $X$-eliminable, then 
$$
  \cov(X\times Y,d_p) \geq \cov(X)\,\cov(Y).
$$
\end{proposition}

\begin{proof}
We argue by induction on $k$. For $k=0$ there is nothing to prove. 
Let $k\geq 1$, write $(B,b):=(B_0,b_0)$ and $A':=A_1=A\setminus B$, 
which has the $Z$-eliminating chain $A_1\supsetneq\cdots\supsetneq A_k$
of length $k-1$. 
Consider a covering of $A\times Y$ by $d_p$-balls
$D_i\subset B_X(x_i)\times B_Y(y_i)$ with centers $(x_i,y_i)\in Z\times Y$, 
$i=1,\dots,N$. Let $I$ be the set of indices $i$ for which $D_i$ meets 
$\{b\}\times Y$. Since $b\in A$, the balls $D_i$ with $i\in I$ cover 
$\{b\}\times Y$. This implies that $b\in B_X(x_i)$ for all
$i\in I$ and $Y\subset \bigcup_{i\in I}B_Y(y_i)$, hence
$|I|\geq\cov(Y)$. 
For each $i\in I$ we have $x_i\in Z$ and $b\in B_X(x_i)$, so the
definition of a $Z$-relative extremal pair implies 
$B_X(x_i)\cap A\subset B$. Thus $D_i\cap(A'\times Y)=\varnothing$, and
therefore the balls $D_i$ with $i\notin I$ cover $A'\times Y$. Their
centers lie in $Z\times Y$, so the induction hypothesis yields
$N-|I| \geq (k-1)\cov(Y)$. Combined with $|I|\geq\cov(Y)$ this yields
$N\geq k\,\cov(Y)$. The remaining assertions follow from
$\cov_Z(A) = k$ and the case $Z=X$. 
\end{proof}

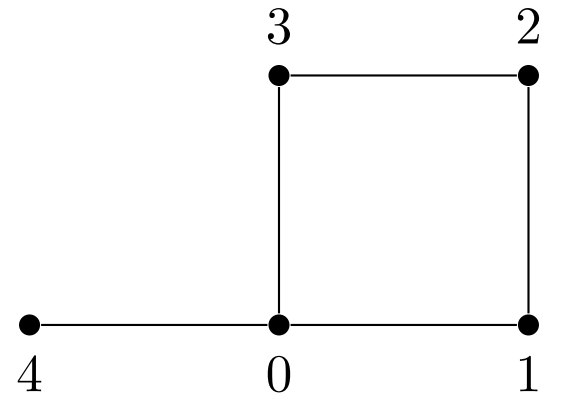
\begin{figure}[ht]
\centering
\begin{tikzpicture}[scale=1.1,every node/.style={font=\small},
    vtx/.style={circle,fill,inner sep=1.4pt}]
  \node[vtx,label=below:$0$] (v0) at (0,0) {};
  \node[vtx,label=below:$1$] (v1) at (1.5,0) {};
  \node[vtx,label=above:$2$] (v2) at (1.5,1.5) {};
  \node[vtx,label=above:$3$] (v3) at (0,1.5) {};
  \node[vtx,label=below:$4$] (v4) at (-1.5,0) {};
  \draw (v0)--(v1)--(v2)--(v3)--(v0)--(v4);
\end{tikzpicture}
\caption{The graph in Example~\ref{ex:eliminable-needed}.}
\label{fig:eliminable-needed}
\end{figure}

\begin{example}\label{ex:eliminable-needed}
Let $X=\{0,1,2,3,4\}$ be the set of vertices of the graph shown in
Figure~\ref{fig:eliminable-needed} with the graph metric of edge length $1$. Let $r=1$,
$p=\infty$, and $Y=A=\{0,1,2,3\}$ with the induced metric.  
Then $(B_X(0),4)$ is an extremal pair in $X$ (so $X$ is
$X$-eliminable) and $\cov(A;X)=\cov(Y)=2$, but the three
$d_\infty$-balls with centers $(1,1)$, $(3,3)$, $(0,0)$ cover $A\times Y$
(see~Example~\ref{ex:Z4} below), 
so $\cov(A\times Y;X\times Y,d_\infty)\leq 3<4=\cov(A;X)\,\cov(Y)$.
Here $A$ has no relative extremal pair, so it is not $X$-eliminable.
\end{example}

{\bf Finite center sets. }
Since infinite sets are often not eliminable
(see Remark~\ref{rem:tree-extremal-false} below), we will apply
Proposition~\ref{prop:extremal} only to finite sets and finitely many
centers. The next lemma ensures the existence of suitable finite sets.

\begin{lemma}\label{lem:finite-witness}
Let $F\subset X$ be finite and $A\subset X$ with
$\cov_F(A)<\infty$. Then there exists a finite subset $A_0\subset A$
with $\cov_F(A_0)=\cov_F(A)$.
\end{lemma}

\begin{proof}
For each subset $G\subset F$ such that the balls $B_X(g)$, $g\in G$, do
not cover $A$, pick a point $a_G\in A\setminus\bigcup_{g\in G}B_X(g)$,
and let $A_0$ be the finite set of these points. Then
$\cov_F(A_0)\leq\cov_F(A)$ because $A_0\subset A$. Conversely, let $G\subset F$ be the set of
centers of a minimal covering of $A_0$ by balls with centers in
$F$. Then the balls with centers in $G$ cover $A$, since otherwise
$a_G\in A_0$ would not be covered. Hence $\cov_F(A)\leq|G| = \cov_F(A_0)$. 
\end{proof}

\begin{definition}\label{def:EP}
A metric space $X$ has {\em property (EP)}\footnote{
``EP'' stands for ``extremal pair''.}
(at radius $r$)
if for all nonempty finite subsets $A,F\subset X$ with 
$A\subset\bigcup_{f\in F}B_X(f)$ the set $A$ has an $F$-relative
extremal pair. 
\end{definition}

\begin{proposition}\label{prop:EP}
Let $X$ be a metric space with property (EP). Then for every
$A\subset X$, every metric space $Y$ and every $p\in[1,\infty]$ we have
$$
  \cov(A\times Y;X\times Y,d_p) \geq \cov(A;X)\,\cov(Y).
$$
In particular, 
$$
  \cov(X\times Y,d_p) \geq \cov(X)\,\cov(Y)\quad\text{and}\quad
  \cov(X\times Y,d_\infty) = \cov(X)\,\cov(Y).
$$
\end{proposition}

\begin{proof}
If $A$ or $Y$ is empty the right hand side is $0$, and if the left
hand side is infinite there is nothing to prove. So assume
$A,Y\neq\varnothing$ and consider a covering of $A\times Y$ by
$N<\infty$ $d_p$-balls with centers $(x_i,y_i)\in X\times Y$. Set
$F:=\{x_1,\dots,x_N\}$. Fix $y\in Y$. For each $a\in A$ the point
$(a,y)$ lies in one of the balls, hence $a\in B_X(x_i)$ for some $i$. 
So $A$ is covered by balls with centers in $F$, and since $F\subset X$ 
we get $\cov(A;X)\leq\cov_F(A)\leq N<\infty$. 
By Lemma~\ref{lem:finite-witness} there exists a finite subset
$A_0\subset A$ with $\cov_F(A_0)=\cov_F(A)$. 
We construct an $F$-eliminating chain for $A_0$: As long as
$A_i\neq\varnothing$, the finite set $A_i\subset A_0$ is covered by
balls with centers in $F$, so by property (EP) it has an $F$-relative
extremal pair $(B_i,b_i)$, and we set $A_{i+1}:=A_i\setminus
B_i$. Since $b_i\in A_i\cap B_i$, the cardinality decreases strictly,
so we reach $A_k=\varnothing$ after finitely many steps. The $N$ given
$d_p$-balls above have centers in $F\times Y$ and cover $A_0\times Y$, so
Proposition~\ref{prop:extremal} with $Z=F$ yields
$$
  N\geq \cov_F(A_0)\,\cov(Y) = \cov_F(A)\,\cov(Y) \geq \cov(A;X)\,\cov(Y)
$$
(and no such covering exists if $\cov(Y)=\infty$ because
$\cov_F(A_0)\geq 1$). Since the covering was arbitrary, this proves the
inequality. For $A=X$ the left hand side equals $\cov(X\times Y,d_p)$, and 
for $p=\infty$ the reverse inequality is Lemma~\ref{lem:basic}(a).
\end{proof}

For the first application, consider nonempty finite subsets $A,F\subset\R$
such that every point of $A$ has distance at
most $r$ from $F$. The {\em leftmost ball $B$ with center in $F$
containing $\max A$} is the closed $r$-ball with center
$c:=\min(F\cap[\max A-r,\infty))$. Then $(B,\max A)$ is an $F$-relative
extremal pair for $A$: if $f\in F$ and $a\in A$ satisfy $|f-\max
A|\leq r$ and $|a-f|\leq r$, then $a\geq f-r\geq c-r$ by minimality of
$c$, and $a\leq\max A\leq c+r$. Hence every subset $X\subset\R$ with
the induced metric has property (EP). 
As a consequence of Proposition~\ref{prop:EP}, we obtain
  
\begin{corollary}\label{cor:R}
Let $X\subset\R$ with the induced metric and $Y$ be an arbitrary
metric space. Then 
$$
  \cov(X\times Y,d_\infty) = \cov(X)\,\cov(Y).
$$
\end{corollary}


The following example shows that Corollary~\ref{cor:R} fails when
$\R$ is replaced by the circle or $\R^2$. 

\begin{example}\label{ex:Z4}
Consider the set $X=\Z/4\Z\subset\R/4\Z$ with the metric induced from
the Euclidean metric on $\R$. Given $p\in(1,\infty]$, pick a radius
$r\in[2^{1/p},2)$. 
Then $r<2$ implies $\cov(X)=2$. On the other hand, in $(X\times
X,d_p)$ we have
$$
  d_p([0,0],[1,1]) = (1^p+1^p)^{1/p} = 2^{1/p}\leq r.
$$
Thus the $d_p$-ball of radius $r$ with center $[i,j]$ consists of all
$[i',j']$ with $\max\{d(i',i),d(j',j)\}\leq 1$, i.e., it agrees with the
$d_\infty$-ball of radius $1$. Therefore, the balls with centers
$[1,1]$, $[3,3]$ and $[0,0]$ cover $X\times X$ and, since $d_p\geq d_\infty$
and by Corollary~\ref{cor:Milan}, we see that
$$
  \cov(X\times X,d_p,r) = 3 < 4 = \cov(X,r)\,\cov(X,r).
$$
See Figure~\ref{fig:torus}.
\begin{figure} [ht]
\hbox{\hspace{2.8cm}
\includegraphics[width=\textwidth]{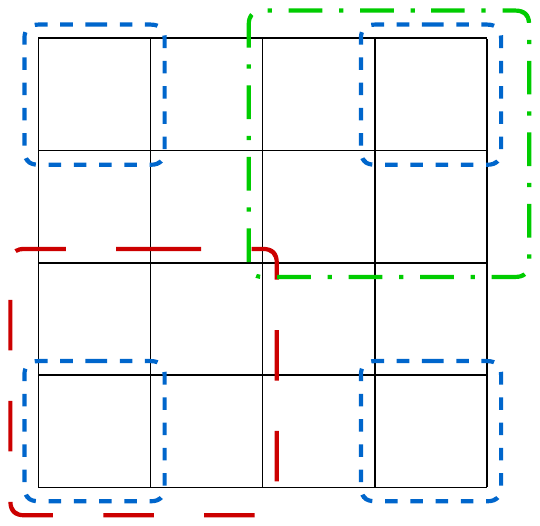}}
\vspace{-12cm}
\caption{Covering $X\times X$ by three balls}
\label{fig:torus}
\end{figure}
Note that this argument works for each $p$ with $1<p\leq\infty$ but
it fails for $p=1$. 
Note also that we can view $X$ as the corners of a square of
side length $1$ with the induced metric from the $\ell^1$-norm on the plane.
\end{example}

For $p=1$, the preceding example is contrasted by the following lemma.

\begin{lemma}\label{lem:cov-two}
Let $X=\Z/4\Z$ be the space from Example~\ref{ex:Z4}. Let $r\in[1,2)$,
so that $\cov(X,r)=2$. Then $\cov(X\times Y,d_1,r)\geq 4$ for each
metric space $Y$ with $\cov(Y,r)\geq 2$. 
\end{lemma}

\begin{proof}
Arguing by contradiction, suppose that $X\times Y$ is covered by $3$
balls of radius $r$. The projection of each ball to $X$ consists of
$3$ points. If two of the projections were equal, then some $x\in X$
would lie in only one projection and the fibre $\{x\}\times Y$ would
be covered by one ball, contradiction $\cov(Y)\geq 2$. Hence, after
cyclic renumberig we may assume that the projections are $\{1,2,3\}$,
$\{2,3,4\}$, and $\{3,4,1\}$, so $X\times Y$ is covered by the
$r$-balls with centers $(2,a)$, $(3,c)$, and $(4,b)$ for some
$a,b,c\in Y$.

Set $\eps:=r-1\in[0,1)$. Looking at the fibres over $1,2,4$ we see that
$$
  Y = B(a,\eps)\cup B(b,\eps)
  = B(a,r)\cup B(c,\eps)
  = B(b,r)\cup B(c,\eps).
$$
If $c\in B(a,\eps)$, then for each $y\in B(c,\eps)$ we have
$$
  d(a,y)\leq d(a,c)+d(c,y)\leq 2\eps<r,
$$
hence $B(c,\eps)\subset B(a,r)$ and the ball $B(a,r)$ covers $Y$,
contradicting $\cov(Y)\geq 2$. Thus $c\notin B(a,\eps)$, and similarly
$c\notin B(b,\eps)$, contradicting $Y = B(a,\eps)\cup B(b,\eps)$. 
\end{proof}

{\bf Subsets of metric trees. }
For a (possibly infinite) graph $G$ we denote by $|G|$ the
associated topological space obtained by gluing intervals
corresponding to the edges at the vertices. A metric on $G$ (i.e., an
assignment of positive lengths to the edges) makes $|G|$ a metric space
with the path metric. 
We begin by recalling the intrinsic characterization of subsets of
metric trees. 

%

Let $X$ be a metric space. Following
Gromov~\cite{Gromov87} (see also~\cite[Ch.~III.H,
Definition~1.19]{Bridson-Haefliger99}), 
we introduce the {\em Gromov product based at $w\in X$}, 
$$
  (x\cdot y)_w := \tfrac12\bigl(d(w,x)+d(w,y)-d(x,y)\bigr). 
$$
The space $X$ is called {\em $0$-hyperbolic} if
\begin{equation}\label{eq:hyp}
  (x\cdot y)_w \geq \min\{(x\cdot z)_w,(y\cdot z)_w\}\qquad\text{for all }x,y,z,w\in X. 
\end{equation}
By~\cite[Corollary 1.1.B]{Gromov87}, condition~\eqref{eq:hyp} holds
for all $w$ if it holds for one $w$. It is equivalent to the 
four-point condition 
$$
  d(x,w)+d(y,z)\leq\max\{d(x,y)+d(z,w),\,d(x,z)+d(y,w)\}
$$ 
for all $x,y,z,w\in X$, see~\cite[Ch.~III.H]{Bridson-Haefliger99} for
its geometric interpretation. 

\begin{lemma}[{\cite[1.1.D, 1.4 and 6.1]{Gromov87}}]\label{lem:Gromov}
Every subset of $|T|$ for a (connected but possibly infinite) metric
tree $T$ is $0$-hyperbolic. A finite metric space is $0$-hyperbolic if
and only if it is isometric to a subset of $|T|$ for a finite metric tree $T$.
\end{lemma}


Let us fix a basepoint $w\in X$ and abbreviate
$$
  |x| := d(w,x),\qquad (x\cdot y) := (x\cdot y)_w =
  \tfrac12\bigl(|x|+|y|-d(x,y)\bigr).   
$$

\begin{lemma}\label{lem:key}
Let $X$ be $0$-hyperbolic and $y,b,c,x\in X$ with $|y|\leq |b|$ and
$|c|\leq |x|$. Then 
$$
  d(y,c)\leq\max\{d(y,x),d(b,c),d(b,x)\}.
$$
\end{lemma}

\begin{proof}
We write $d(y,c)=|y|+|c|-2(y\cdot c)$. If $(y\cdot c)\geq(y\cdot x)$, then 
$d(y,c)\leq |y|+|x|-2(y\cdot x)=d(y,x)$. If $(y\cdot c)\geq(b\cdot c)$, then 
$d(y,c)\leq |b|+|c|-2(b\cdot c)=d(b,c)$. 
In the remaining case $(y\cdot c)<(y\cdot x)$ and $(y\cdot c)<(b\cdot c)$ we claim that
$(x\cdot b)\leq(y\cdot c)$. Otherwise condition~\eqref{eq:hyp} would give
$(x\cdot c)\geq\min\{(x\cdot b),(b\cdot c)\}>(y\cdot c)$ and then
$(y\cdot c)\geq\min\{(y\cdot x),(x\cdot c)\}>(y\cdot c)$, a contradiction. Hence 
$d(y,c)\leq |b|+|x|-2(x\cdot b)=d(b,x)$. 
\end{proof}

\begin{proposition}\label{prop:tree-extremal}
Every $0$-hyperbolic metric space $X$ has
property (EP) of Definition~\ref{def:EP} (at every radius $r>0$).
\end{proposition}

\begin{proof}
Consider nonempty finite subsets $A,F\subset X$ with 
$A\subset\bigcup_{f\in F}B_X(f)$, i.e., each point of $A$ has
distance at most $r$ from $F$. As above, we fix a basepoint $w\in X$
and write $|x|=d(w,x)$. Let $b\in A$
maximize $|\ |$ on $A$, and let $c$ minimize $|\ |$ on the set
$\{f\in F\mid d(b,f)\leq r\}$
(which is nonempty because $b$ has distance at most $r$ from $F$).
It follows that $d(b,c)\leq r$.
We claim that $(B_X(c),b)$ is an $F$-relative extremal pair for $A$. 

To see this, consider $f\in F$ with $b\in B_X(f)$ and $a\in
B_X(f)\cap A$, so $d(b,f)\leq r$ and
$d(a,f)\leq r$. Then $|a|\leq |b|$ by the choice of $b$, and
$|c|\leq |f|$ by the choice of $c$. So Lemma~\ref{lem:key} yields
$$
  d(a,c)\leq\max\{d(a,f),d(b,c),d(b,f)\}\leq r,
$$
thus $a\in B_X(c)$. This proves the claim and therefore the proposition.
\end{proof}

By Proposition~\ref{prop:tree-extremal}, Lemma~\ref{lem:Gromov} and
Proposition~\ref{prop:EP} we obtain

\begin{corollary}\label{cor:tree-subset}
Let $T$ be a metric tree with arbitrary edge lengths, and
$X\subset |T|$ any subset with the induced metric. Fix a radius $r>0$
and let $Y$ be an arbitrary metric space. Then  
$$
  \cov(X\times Y,d_p) \geq \cov(X\times Y,d_\infty) =
  \cov(X)\,\cov(Y) \quad\text{for all }p\in[1,\infty]. 
$$
In particular, Conjecture~\ref{conj:cov} holds if $X$ is a subset of a
metric tree, and for the vertex set $X$ of a finite tree with $r=1$ it
gives Vizing's conjecture when one of the graphs is a tree.
\end{corollary}

\begin{remark}\label{rem:tree-extremal-false}
Infinite subsets of a tree need not have extremal pairs. Let $|T|=[0,10]$,
$X=(0,10)$ and $r=1$. For $b\in X$ and $y\in X$ with $|y-b|<2$, the
midpoint $m=(y+b)/2$ lies in $X$ and $y,b\in B_X(m)$. So the union of
all balls containing $b$ contains $X\cap(b-2,b+2)$, which has diameter
$>2$. Since every ball $B_X(c)\subset[c-1,c+1]$ has diameter at most $2$,
$X$ has no extremal pair. This also shows that $\bar X=[0,10]$ has
extremal pairs but is not $\bar X$-eliminable: every chain of removals of
relative extremal pairs gets stuck at $(2,8)$. This is why we only use
finite sets and finite center sets in Proposition~\ref{prop:EP}.
\end{remark}

\begin{remark}\label{rem:tree-literature}
For vertex sets of graphs with the graph metric and integer radius
$r=k$, cases of Corollary~\ref{cor:tree-subset} are known in graph
theory. Meir and Moon~\cite{MeirMoon75} proved that for a tree the
minimal size of a set of vertices such that every vertex has distance
at most $k$ from it equals the maximal number of vertices with pairwise
distances $>2k$; their proof uses the vertex at distance $k$ from the
end of a longest path, similarly to Proposition~\ref{prop:tree-extremal}. 
Together with a packing lower bound, this yields
$\gamma_k(G\square H)\geq\gamma_k(G)\,\gamma_k(H)$ for the distance-$k$
domination numbers if $G$ is a tree, and more generally if $G$ is
strongly chordal~\cite{BlidiaChellali07}. For $k=1$ this is the case of
Vizing's conjecture where one graph is a tree, due to Jacobson and
Kinch (\cite{Jacobson-Kinch86}, see also~\cite{Bresar-et-al12}). 
\end{remark}

\section{Covering density}\label{sec:density}

In this section we prove an infinitesimal version of the covering
inequality~\eqref{eq:cov}. See~\cite{FejesToth17,Conway-Sloane99} for background. 

By a {\em convex body} $A\subset \R^m$ we will mean a compact convex subset containing
$0$ in its interior. We associate to $A$ its (nonsymmetric) Minkowski norm
$$
  \|x\|_A := \inf\{t>0\mid t^{-1}x\in A\}. 
$$
It satisfies $\|x\|_A=0\Longleftrightarrow x=0$,
$\|x+y\|_A\leq\|x\|_A+\|y\|_A$, and $\|sx\|_A=s\|x\|_A$ for $s\geq 0$;
it is an honest norm iff $A$ is centrally symmetric, i.e., $A=-A$.
The unit balls of the corresponding (nonsymmetric) Minkowski metric
$d_A(x,y):=\|y-x\|_A$ are the translates of $A$.  

{\bf Volume. }
For a compact subset $X\subset\R^m$ we denote by $\vol(X)$ its
Lebesgue measure. 
Consider now two convex bodies $A\subset\R^m$ and $B\subset\R^n$ with
associated norms $\|\ \|_A$ and $\|\ \|_B$. For $1\leq p<\infty$, the unit
ball of the norm 
$$
  \|(x,y)\|_p = (\|x\|_A^p+\|y\|_B^p)^{1/p}
$$
on $\R^{m+n}$ is
\begin{align*}
  A\times_p B
  &:= \{(x,y)\in\R^{m+n}\;\bigl|\;\|x\|_A^p+\|y\|_B^p\leq 1\} \cr
  &= \{(sa,tb) \mid a\in\p A,\,b\in\p B,\,s,t\geq 0,\,s^p+t^p\leq 1\}. 
\end{align*}
Assume for the moment that $A$ and $B$ have smooth boundaries. Then
$$  
  \vol(sA) = s^m\vol(A),\qquad \frac{d}{ds}\vol(sA)=ms^{m-1}\vol(A),
$$
and similarly for $B$. Thus
\begin{align*}
  \vol(A\times_p B)
  &= \iint_{s,t\geq 0,\,s^p+t^p\leq 1}ms^{m-1}\vol(A)nt^{n-1}\vol(B)ds\,dt \cr
  &= \vol(A)\vol(B)\,I_{m,n,p}
\end{align*}
with
$$
  I_{m,n,p}
  := \int_{s=0}^1ms^{m-1}\Bigl(\int_{t=0}^{(1-s^p)^{1/p}}nt^{n-1}dt\Bigr)ds
  = \int_0^1ms^{m-1}(1-s^p)^{n/p}ds.
$$
By an approximation argument, the relation $\vol(A\times_p B) =
\vol(A)\vol(B)\,I_{m,n,p}$ continues to hold for arbitrary convex
bodies $A,B$. It also holds for $p=\infty$ with $A\times_\infty
B=A\times B$ and $I_{m,n,\infty}=1$. 
Note that the above integral for $I_{m,n,p}$ is defined for all
positive real numbers $m,n,p$.

\begin{lemma}\label{lem:I}
For all $m,n,p>0$ we have
\begin{equation}\label{eq:I}
  I_{m,n,p} = \frac{\Gamma\Bigl(\frac{m}{p}+1\Bigr)
    \Gamma\Bigl(\frac{n}{p}+1\Bigr)}{\Gamma\Bigl(\frac{m+n}{p}+1\Bigr)}.  
\end{equation}
Moreover, this quantity satisfies $\frac{\p}{\p n}I_{m,n,p}<0$ for all $m,n,p>0$.
\end{lemma}
  
\begin{proof}
Formula~\eqref{eq:I} follows by the following computation:
\begin{align*}
  I_{m,n,p} &= \frac{m}{p}\int_0^1t^{\frac mp-1}(1-t)^{\frac np}dt 
  = \frac{m}{p} B\Bigl(\frac{m}{p},\frac{n}{p}+1\Bigr) \cr
  &= \frac{m}{p}\,\frac{\Gamma\bigl(\frac mp\bigr) \Gamma\bigl(\frac
    np+1\bigr)}{\Gamma\bigl(\frac{m+n}{p}+1\bigr)}   
  = \frac{\Gamma\bigl(\frac mp+1\bigr)\Gamma\bigl(\frac
    np+1\bigr)}{\Gamma\bigl(\frac{m+n}{p}+1\bigr)}. 
\end{align*}
Here the first equality is obtained by substituting $t=s^p$; 
the second one by definition of Euler's beta function;
and the third one by the relation between the beta and the gamma function.

For the last assertion, consider for fixed $y>0$ the function of $x>0$,
$$
  f(x)
  := \log\Bigl(\frac{\Gamma(x+y+1)}{\Gamma(x+1)\Gamma(y+1)}\Bigr)
  = \log\Gamma(x+y+1)-\log\Gamma(x+1)-\log\Gamma(y+1).
$$
Then $f'(x) = \psi(x+y+1)-\psi(x+1) > 0$ because $y>0$ and the digamma
function $\psi(x):=\frac{d}{dx}\log\Gamma(x)$ is strictly increasing.
\end{proof}

{\bf Covering density. }
For a convex body $A\subset\R^m$ and a compact subset $X\subset\R^m$ let
$$
  \vartheta_T(X;A) := \inf_\mathcal{C}\frac{N\,\vol(A)}{\vol(X)},
$$
where the infimum is taken over all coverings $\mathcal{C}$ of $X$ by
translates $c_1+A,\dots,c_N+A$, $c_i\in\R^m$. Denote by
$B^m(R)\subset\R^m$ the closed Euclidean ball of radius $R$ around the
origin and write $B^m=B^m(1)$.
The {\em translational covering density} of $A$ is 
$$
  \vartheta_T(A) := \liminf_{R\to\infty}\vartheta_T(B^m(R);A)
  = \liminf_{r\to 0}\vartheta_T(B^m(1);rA).
$$

\begin{proposition}[Groemer~\cite{Groemer63}]\label{prop:Groemer}
For every compact subset $X\subset\R^m$ with
\begin{equation}\label{eq:X}
  \vol(X)=\vol(\Int X)>0
\end{equation}
we have
$$
  \vartheta_T(A) = \lim_{r\to 0}\vartheta_T(X;rA).
$$
\end{proposition}
  
\begin{remark}
For a compact subset $X\subset\R^m$ we have the equivalences
$$
  \vol(X)=\vol(\Int X) \Longleftrightarrow \vol(\p X)=0
  \Longleftrightarrow \text{$X$ is Jordan measurable},
$$
and these conditions are satisfied if $X$ is convex
(\cite{Minkowski03}, see also~\cite{Szabo97}).
\end{remark}

Analogously one defines the {\em lattice covering density}
$\vartheta_L(A)$ where one only allows translations in some lattice,
and the {\em covering density} $\vartheta(A)$ where one allows all
translations as well as rotations. Thus 
$$
  \vartheta(A) \leq \vartheta_T(A) \leq \vartheta_L(A).
$$
We are primarily interested in $\vartheta_T$ because by definition it
is related to the covering number $\cov^{\rm ext}(X,d_A,r)$, the
minimal number of balls $c+rA$ with centers $c\in\R^m$ anywhere
covering $X$, by
\begin{equation}\label{eq:cov-number-density}
  \cov^{\rm ext}(X,d_A,1) = \frac{\vol(X)}{\vol(A)}\vartheta_T(X;A).
\end{equation}
Since $\cov(X,d_A,1)$ only allows balls centered in $X$, we have
$\cov(X,d_A,1)\geq\cov^{\rm ext}(X,d_A,1)$.
On the other hand, we have

\begin{lemma}\label{lem:cov-ext}
If $X$ satisfies~\eqref{eq:X}, then
$$
  \frac{\cov(X,d_A,r)}{\cov^{\rm ext}(X,d_A,r)}\to 1\quad\text{as } r\to 0.
$$
\end{lemma}

\begin{proof}
Choose $0<\eps\leq R$ with $\eps B^m\subset A\subset RB^m$,
and take an optimal covering of $X$ by balls $c_k+rA$, $1\leq k\leq N$, each
meeting $X$. If $c_k\notin X$, then the segment from $c_k$ to a point of
$X\cap(c_k+rA)$ lies in $c_k+rA$ and meets $\p X$, so $c_k+rA$ lies in
the Euclidean $2Rr$-neighbourhood $(\p X)_{2Rr}$ of $\p X$. Hence $X$
is covered by the balls with $c_k\in X$ together with
$Z_r:=X\cap(\p X)_{2Rr}$. A maximal subset $P\subset Z_r$ with
pairwise Euclidean distances $\geq\eps r$ satisfies
$Z_r \subset \bigcup_{p\in P}(p+\eps rB^m) \subset \bigcup_{p\in
  P}(p+rA)$ with centers $p\in X$, and comparing volumes of the
disjoint balls $p+\frac{\eps r}{2}\Int B^m$ gives
$$
  |P|\leq\frac{2^m\vol\bigl((\p X)_{(2R+\eps)r}\bigr)}{\eps^mr^m\vol(B^m)}.
$$
Since $N\geq\vol(X)/(r^m\vol(A))$, we get
\begin{align*}
  1\leq\frac{\cov(X,d_A,r)}{\cov^{\rm ext}(X,d_A,r)}
  &\leq 1+\frac{|P|}{N}\cr
  &\leq 1+\frac{2^m\vol(A)}{\eps^m\vol(B^m)\vol(X)}\,
  \vol\bigl((\p X)_{(2R+\eps)r}\bigr)\xrightarrow[r\to 0]{} 1,
\end{align*}
since $\p X$ is compact and $\vol(\p X)=0$ by~\eqref{eq:X}.
\end{proof}

Similarly, for a convex body $B\subset\R^n$ and a compact subset
$Y\subset\R^n$ we have
\begin{equation*}
  \cov^{\rm ext}(Y,d_B,1) = \frac{\vol(Y)}{\vol(B)}\vartheta_T(Y;B).
\end{equation*}
For $p\in[1,\infty]$ it follows that
\begin{align*}
  \cov^{\rm ext}(X\times Y,d_p,1)
  &= \frac{\vol(X\times
    Y)}{\vol(A\times_pB)}\vartheta_T(X\times Y;A\times_pB) \cr
  &= \frac{\vol(X)\vol(Y)\vartheta_T(X\times Y;A\times_pB)}{\vol(A)\vol(B)\,I_{m,n,p}},
\end{align*}
Rescaling with $r>0$ and using Proposition~\ref{prop:Groemer}, this yields
\begin{align*}
  \frac{\cov^{\rm ext}(X\times Y,d_p,r)}{\cov^{\rm ext}(X,d_A,r)\cov^{\rm ext}(Y,d_B,r)}
  &= \frac{\vartheta_T(X\times Y;rA\times_prB)}{I_{m,n,p}\vartheta_T(X;rA)\vartheta_T(Y;rB)}\cr
  &\xrightarrow[r\to 0]{} \frac{\vartheta_T(A\times_pB)}{I_{m,n,p}\vartheta_T(A)\vartheta_T(B)}
\end{align*}
provided that $X,Y$ satisfy~\eqref{eq:X}. Then $X\times Y$
satisfies~\eqref{eq:X} as well, so by Lemma~\ref{lem:cov-ext} (applied
to $X$, $Y$, and $X\times Y$ with the body $A\times_pB$) the same limit
holds for the ratio with $\cov$ in place of $\cov^{\rm ext}$. So we have shown

\begin{lemma}\label{lem:cov-infinitesimal}
Suppose that convex bodies $A\subset\R^m$, $B\subset\R^n$ satisfy the
{\em infinitesimal covering inequality}  
\begin{equation}\label{eq:cov-infinitesimal}
  \frac{\vartheta_T(A\times_pB)}{I_{m,n,p}} > \vartheta_T(A)\vartheta_T(B).
\end{equation}
Then for all compact subsets $X\subset\R^m$, $Y\subset\R^n$
satisfying~\eqref{eq:X} there exists $r_{X,Y}>0$ such that for all
$0<r<r_{X,Y}$ the covering inequality~\eqref{eq:cov} holds
with $d_A,d_B$ in place of $d_X,d_Y$. \qed
\end{lemma}

Using Lemma~\ref{lem:I}, we can rewrite inequality~\eqref{eq:cov-infinitesimal} as
\begin{equation}\label{eq:cov-infinitesimal-Gamma}
  \frac{\Gamma\Bigl(\frac{m+n}{p}+1\Bigr)}{\Gamma\Bigl(\frac{m}{p}+1\Bigr)
  \Gamma\Bigl(\frac{n}{p}+1\Bigr)}\vartheta_T(A\times_pB) > \vartheta_T(A)\vartheta_T(B).
\end{equation}

For $p=1$ this leads to the following infinitesimal version of
Conjecture~\ref{conj:cov}. 

\begin{conjecture}\label{conj:cov-infinitesimal}  
For all convex bodies $A\subset\R^m,B\subset\R^n$ we have
\begin{equation}\label{eq:cov-infinitesimal-binom}
  \binom{m+n}{n}\vartheta_T(A\times_1B) > \vartheta_T(A)\vartheta_T(B).
\end{equation}
\end{conjecture}

\begin{remark}\label{rem:cov-infinitesimal-1D}
The non-strict version ($\geq$ instead of $>$) of
inequality~\eqref{eq:cov-infinitesimal} holds for $n=1$ and
each $p\in[1,\infty]$ because then $\vartheta_T(B)=1$ and 
$$
  \frac{\vartheta_T(A\times_pB)}{I_{m,n,p}} \geq \vartheta_T(A\times B)
  = \vartheta_T(A).
$$
Here the inequality holds because $A\times_pB\subset A\times B$
and $\vol(A\times_pB)=I_{m,1,p}\vol(A\times B)$, so every covering by
translates of $A\times_pB$ is a covering by the same number of
translates of $A\times B$.
For the equality, ``$\leq$'' follows by covering a closed interval $X\subset\R$ by
translates of the closed interval $rB$ with disjoint interiors, and
``$\geq$'' by slicing: if $N$ translates of $r(A\times B)$ cover $Z\times J$
for a compact convex $Z\subset\R^m$ and an interval $J$, then for each $t\in J$
the translates whose second factor contains $t$ yield a covering of $Z$ by
translates of $rA$; integrating their number over $t\in J$ gives
$N\,\vol(rB)\geq\vol(J)\vol(Z)\vartheta_T(Z;rA)/\vol(rA)$, i.e.,
$\vartheta_T(Z\times J;r(A\times B))\geq\vartheta_T(Z;rA)$, and
Proposition~\ref{prop:Groemer} yields $\vartheta_T(A\times B)\geq\vartheta_T(A)$.
Since $\vartheta_T(A\times_pB)\geq 1$, the strict
inequality~\eqref{eq:cov-infinitesimal} holds for $n=1$ whenever
$\vartheta_T(A)<1/I_{m,1,p}$. For $p=1$ we have $1/I_{m,1,1}=m+1$, so
Conjecture~\ref{conj:cov-infinitesimal} holds for $n=1$ and $m\leq 2$:
for $m=1$ we have $\vartheta_T(A)=1<2$, and for $m=2$ F\'ary's bound
from~\cite{Fary50} yields $\vartheta_T(A)\leq 3/2<3$.
\end{remark}

\begin{remark}\label{rem:Zong}
For $p=\infty$, the non-strict version of~\eqref{eq:cov-infinitesimal}
becomes 
$$
  \vartheta_T(A\times B) \geq \vartheta_T(A)\vartheta_T(B).
$$
C.~Zong~\cite{Zong05} has shown that for every
convex body $A$ with $\vartheta_T(A)>1$ there exists a convex body $B$
such that $\vartheta_T(A\times B) < \vartheta_T(A)\vartheta_T(B)$. 
Hence, the non-strict version of~\eqref{eq:cov-infinitesimal} fails
in general for sufficiently large $p$. 
\end{remark}

{\bf Estimates on the covering density. }
We will use the trivial lower bound $\vartheta_T(A\times_1B)\geq 1$ and
replace~\eqref{eq:cov-infinitesimal-binom} by the stronger inequality
\begin{equation}\label{eq:cov-infinitesimal-binom-strong}
  \binom{m+n}{n} > \vartheta_T(A)\vartheta_T(B).
\end{equation}
On the other hand, C.\,Rogers~\cite{Rogers57} has proved a general upper bound
for each convex body $A\subset\R^m$, $m\geq 3$, where $\log$ denotes
the natural logarithm: 
\begin{equation}\label{eq:rogers}
\begin{aligned}
  \vartheta_T(A)
  &\leq \Bigl(1+\frac{1}{m\log m}\Bigr)^m\Bigl(1+m\log\bigl(m\log m\bigr)\Bigr) \cr
  &\leq m\log m + m\log\log m + 5m.
\end{aligned}
\end{equation}
Using for $m=1,2$ the bounds from Remark~\ref{rem:cov-infinitesimal-1D}, 
we define the functions
\begin{align*}
\wt\rho(m) &:= \begin{cases}
  1 & m=1, \cr
  3/2 & m=2, \cr
  \Bigl(1+\frac{1}{m\log m}\Bigr)^m\Bigl(1+m\log\bigl(m\log m\bigr)\Bigr) & m\geq 3,
\end{cases} \cr
\rho(m) &:= \begin{cases}
  1 & m=1, \cr
  3/2 & m=2, \cr
  m\log m + m\log\log m + 5m & m\geq 3.
\end{cases}
\end{align*}
Then $\wt\rho(m)\leq\rho(m)$, and~\eqref{eq:cov-infinitesimal-binom}
is implied by 
\begin{equation}\label{eq:cov-infinitesimal-binom-rogers-improved}
  \binom{m+n}{n} > \wt\rho(m)\wt\rho(n),
\end{equation}
which is in turn implied by
\begin{equation}\label{eq:cov-infinitesimal-binom-rogers}
  \binom{m+n}{n} > \rho(m)\rho(n).
\end{equation}

\begin{proposition}\label{prop:cov-infinitesimal}
Conjecture~\ref{conj:cov-infinitesimal} is true if $\min\{m,n\}>5$.
More precisely, it holds for all pairs $(n,m)$ with $1\leq n\leq m$ 
except possibly the following ones: $(1,m)$ with $m\geq 3$,
$(2,m)$ with $3\leq m\leq 12$, $(3,m)$ with
$3\leq m\leq 15$, $(4,m)$ with $4\leq m\leq 8$, $(5,5)$, and $(5,6)$.
\end{proposition}

\begin{proof}
We will prove that
inequality~\eqref{eq:cov-infinitesimal-binom-rogers-improved} holds
for the pairs in the proposition. In view of
Remark~\ref{rem:cov-infinitesimal-1D} we may assume $2\leq n\leq m$.
For $m\geq n\geq 4$ we estimate $\rho(m)^2\geq\rho(m)\rho(n)$ and 
$$
  \binom{m+n}{n}
  = \frac{m+n}{n} \cdot\frac{m+n-1}{n-1} \cdots \frac{m+1}{1}
  \geq \binom{m+4}{4} \geq \frac{m^4}{24}.
$$
Thus~\eqref{eq:cov-infinitesimal-binom-rogers} is implied by
$$
  \sigma(m) := \frac{\sqrt{24}\rho(m)}{m^2} = \frac{\sqrt{24}(\log m +
    \log\log m + 5)}{m} < 1.
$$
The function $\sigma$ decreases and $\sigma(51)\approx 0.99<1$,
so~\eqref{eq:cov-infinitesimal-binom-rogers} holds if $n\geq 4$ and
$m>50$.
One directly checks that~\eqref{eq:cov-infinitesimal-binom-rogers}
also holds for $n\in\{2,3\}$ and $m>50$. In the remaining cases $2\leq 
n\leq m\leq 50$, one checks
that~\eqref{eq:cov-infinitesimal-binom-rogers-improved} holds except for the
cases listed in the proposition. 
\end{proof}

The proofs actually give the following more quantitative version of
Proposition~\ref{prop:cov-infinitesimal} and
Lemma~\ref{lem:cov-infinitesimal}.

\begin{corollary}\label{cor:cov-infinitesimal}
For $(n,m)$ as  
in Proposition~\ref{prop:cov-infinitesimal} there exist constants
$\mu_{m,n}>1$ and $p_{m,n}>1$ depending only on
$m,n$ such that for all convex bodies $A\subset\R^m$, $B\subset\R^n$
and $1\leq p\leq p_{m,n}$ we have
\begin{equation*}
  \frac{\vartheta_T(A\times_pB)}{I_{m,n,p}} \geq \mu_{m,n}^4\vartheta_T(A)\vartheta_T(B).
\end{equation*}
Given such $p,A,B$, for all compact subsets $X\subset\R^m$, $Y\subset\R^n$
satisfying~\eqref{eq:X} there exists $r_{X,Y}>0$ such that for all
$0<r\leq r_{X,Y}$ we have
\begin{equation*}
  \cov(X\times Y,d_p,r) \geq \mu_{m,n}^3\cov(X,d_A,r)\,\cov(Y,d_B,r).
\end{equation*}
\end{corollary}

To deal with the remaining cases listed in
Proposition~\ref{prop:cov-infinitesimal} one needs a better
lower bound on $\vartheta_T(A\times_1B)$, or better upper bounds on
$\vartheta_T(A)$ in the dimensions $3\leq m\leq 15$. Here is
what appears to be known about this, cf.~\cite{FejesToth22}.

$\bullet$ For {\em centrally symmetric} $A\subset\R^3$, E.~Smith~\cite{Smith00} has
proved the upper bound $\vartheta_T(A)\leq\vartheta_L(A)\leq
4\delta_L(A)\leq 4$ (where $\delta_L(A)\leq 1$ is the lattice packing
density).\footnote{
See R.\,Prosanov~\cite{Prosanov} for a kind of duality between
the translational packing density $\delta_T$ and the translational
covering density $\vartheta_T$ in the case that both are close to $1$.
}
Replacing $\wt\rho(3)\approx 10.14$ by $4$ shows that, in the
centrally symmetric case, Conjecture~\ref{conj:cov-infinitesimal} also
holds in the cases $(2,3)$ and $(3,m)$ with $8\leq m\leq 15$.

$\bullet$ P.\,Erd\H os and C.\,Rogers~\cite{Erdos-Rogers62} have
improved Rogers' weaker bound for large $m$ to
\begin{align*}
  \vartheta_T(A) \leq m\log m + m\log\log m + 4m,
\end{align*}
but this doesn't help with the cases $m\leq 15$. 

$\bullet$ G.~Fejes T\'oth~\cite{FejesToth09} has proved for convex
$A\subset\R^m$ and large $m$ the bound
\begin{align*}
  \vartheta_T(A)
  &\leq \Bigl(1+\frac{1}{m\log m}\Bigr)^mm\Bigl(\log m+\log\log
  m+\frac{1}{\log m}\Bigr).
\end{align*}
This is better than~\eqref{eq:rogers} for large $m$, but for $m\leq
15$ (assuming it holds there as well) it would be worse
than~\eqref{eq:rogers}. 

$\bullet$ J.\,Rolfes and F.\,Vallentin~\cite{Rolfes-Vallentin18} present an
interesting general approach to estimating covering numbers. They
consider a compact metric space $(X,d)$ equipped with a Borel measure $\om$
satisfying $\om(X)=1$ and the following conditions on closed balls:\\
\ \ \ (a) $\om(B(x,r))=\om(B(y,r))$ for all $x,y\in X$ and $r\geq 0$;\\
\ \ \ (b) $\om(B(x,r))>0$ for all $x\in X$ and $r>0$. \\
Denoting $\om_r:=\om(B(x,r))$, it is proved in~\cite{Rolfes-Vallentin18}
that for each $r/2>\eps>0$ we have
\begin{equation}\label{eq:Rolfes-Vallentin18}
  \frac{1}{\om_r}\leq \cov(X,d,r) \leq
  \frac{1}{\om_{r-\eps}}
  \Bigl(\log\Bigl(\frac{\om_{r-\eps}}{\om_\eps}\Bigr)+1\Bigr).  
\end{equation}
We apply this to a ball $X\subset\R^m$ with the Euclidean volume
$\om$ such that $\om(X)=1$, and $d=d_A$ for a centrally symmetric
convex body $A\subset\R^m$. For small $r$ there is a small error in
condition (a) that vanishes as $r\to 0$. Then $\om_r=r^m\vol(A)$
and~\eqref{eq:Rolfes-Vallentin18} becomes
\begin{equation*}
  \frac{1}{r^m\vol(A)}\leq \cov(X,d_A,r) \leq
  \frac{1}{(r-\eps)^m\vol(A)}
  \Bigl(\log\Bigl(\frac{(r-\eps)^m}{\eps^m}\Bigr)+1\Bigr).  
\end{equation*}
Combining this with the relation
\begin{equation*}
  \cov^{\rm ext}(X,d_A,r) = \frac{\vol(X)}{r^m\vol(A)}\vartheta_T(X;rA)
\end{equation*}
from~\eqref{eq:cov-number-density} (note that $\cov^{\rm ext}\leq\cov$, and the lower bound
also holds for $\cov^{\rm ext}$ by comparing volumes) and writing $\eps=xr$ with
$0<x<1/2$, we obtain
\begin{align*}
  1 \leq \vartheta_T(X;rA)
  &\leq \Bigl(\frac{r}{r-\eps}\Bigr)^m
  \Bigl(\log\Bigl(\frac{(r-\eps)^m}{\eps^m}\Bigr)+1\Bigr) \cr
  &= \frac{1}{(1-x)^m}
  \Bigl(m\log\Bigl(\frac{1}{x}-1\Bigr)+1\Bigr) =: f_m(x).
\end{align*}
In the limit $r\to 0$ this yields for centrally symmetric $A$ the
upper bound 
\begin{equation}\label{eq:Rolfes-Vallentin-density}
  \vartheta_T(A) \leq \inf_{0<x<1/2}f_m(x).
\end{equation}
Numerical evaluation shows that at least for small $m$ this bound is
slightly better than~\eqref{eq:rogers}, but it does not lead to an
improvement of Proposition~\ref{prop:cov-infinitesimal}. 
M.\,Nasz\'odi~\cite{Naszodi16} derives a slightly weaker bound
than~\eqref{eq:Rolfes-Vallentin-density} that also holds for
nonsymmetric $A$.

$\bullet$ Since for each $m\in\N$ there exist convex bodies $C\subset\R^m$ with
$\vartheta_T(C)=1$ (e.g.~cubes), there is no general lower bound on
the covering number beyond $\vartheta_T\geq 1$. However, the following
question might still have a positive answer.

\begin{question}
Is there a nontrivial lower bound on $\vartheta_T(A\times_1B)$ in
terms of $\vartheta_T(A)$ and $\vartheta_T(B)$?
\end{question}

{\bf Ellipsoids. }
Let us now consider the special case that $A,B$ are ellipsoids. The
following lemma reduces this to the case of balls.

\begin{lemma}\label{lem:ellipsoid-to-ball}
For each convex body $A\subset\R^m$ and linear isomorphism
$\phi:\R^m\to\R^m$ we have $\vartheta_T(\phi(A))=\vartheta_T(A)$.
\end{lemma}

\begin{proof}
If a compact set $X\subset\R^m$ is covered by translates
$c_1+A,\dots,c_N+A$, then $\phi(X)$ is covered by translates
$\phi(c_1)+\phi(A),\dots,\phi(c_N)+\phi(A)$ and vice versa, hence
$\vartheta_T(\phi(X);\phi(A))=\vartheta_T(X;A)$. The statement now
follows with $X=B^m(R)$ in the limit $R\to\infty$: with the ellipsoid
$E:=\phi^{-1}(B^m(1))$ we get
$\vartheta_T(B^m(R);\phi(A))=\vartheta_T(RE;A)=\vartheta_T(E;R^{-1}A)$, which tends
to $\vartheta_T(A)$ by Proposition~\ref{prop:Groemer}.
\end{proof}

For balls there are much better upper bounds on $\vartheta_L$ and thus
on $\vartheta_T$,
see~\cite[\S1.3]{Schuermann09} for the best known bounds in dimensions $1$ to $24$.
With these bounds the remaining cases in
Proposition~\ref{prop:cov-infinitesimal} hold and we obtain

\begin{corollary}\label{cor:cov-infinitesimal-ellipsoid}
For ellipsoids $A\subset\R^m$ and $B\subset\R^n$,
Conjecture~\ref{conj:cov-infinitesimal} is true for all $m\geq n\geq 2$.
For $n=1$ it is true if $m\leq 17$, $m=20$, $m=23$ or $m=24$, since
then the bounds from~\cite[\S1.3]{Schuermann09} 
give $\vartheta_T(B^m)<m+1$ (see Remark~\ref{rem:cov-infinitesimal-1D}).
\qed
\end{corollary}


The $2$-product $A\times_2 B$ of ellipsoids $A,B$ is again an
ellipsoid. This raises the question
whether~\eqref{eq:cov-infinitesimal} holds for ellipsoids
$A\subset\R^m$ and $B\subset\R^n$ with $p=2$, 
\begin{equation}\label{eq:cov-infinitesimal-ellipsoid}
  \frac{\vartheta_T(A\times_2B)}{I_{m,n,2}} > \vartheta_T(A)\vartheta_T(B).
\end{equation}
Lemma~\ref{lem:ellipsoid-to-ball} reduces the question to balls
$B^d\subset\R^d$, so~\eqref{eq:cov-infinitesimal-ellipsoid} is
equivalent to  
\begin{equation}\label{eq:cov-infinitesimal-ball}
  \frac{\vartheta_T(B^{m+n})}{I_{m,n,2}} > \vartheta_T(B^m)\vartheta_T(B^n).
\end{equation}
Computations using the values for $\vartheta_T(B^d)$
from~\cite[\S1.3]{Schuermann09} show that~\eqref{eq:cov-infinitesimal-ball}  
fails in certain dimensions (e.g.~$m=21, n=3$) if we replace
$\vartheta_T(B^{m+n})$ by its trivial lower bound $1$. 
On the other hand, Coxeter, Few and Rogers~\cite{Coxeter-Few-Rogers59}
have proved the lower bound 
$$ 
  \vartheta_T(B^n) \geq \tau_n = \Bigl(\frac{2n}{n+1}\Bigr)^{n/2}\sigma_n.
$$
Here for $\tau_n$ one considers a regular simplex in $\R^n$
of edge length $\sqrt{2(n+1)/n}$ and its covering by the $n+1$ balls of
radius $1$ around its vertices; then $\tau_n\geq 1$ is the sum of the 
volumes of the intersections of the balls with the simplex divided by
the volume of the simplex.  
For $\sigma_n$ one considers a regular simplex $\Delta\subset\R^n$
of edge length $2$ and its packing by the $n+1$ balls of
radius $1$ around its vertices; then $\sigma_n\leq 1$ is the sum of the 
volumes of the intersections of the balls with the simplex divided by
the volume of the simplex.  

The quantity $\sigma_n$ is defined and estimated by Rogers~\cite{Rogers58}.
By formulas (5) and (6) of that paper we have
$$
  \sigma_n = \frac{(n+1)!}{2^{n/2}}V_n,
$$
where $V_n$ is the volume of the region in $\R^n$ defined by the inequalities
$$
  \sum_{j=1}^nx_j^2 + \Bigl(\sum_{j=1}^nx_j\Bigr)^2\leq 1, \quad
  x_j\geq 0\text{ for }1\leq j\leq n.
$$
Equivalently, $V_n$ can be computed as the integral
$$
  V_n = \frac{1}{\Gamma(n/2+1)}\int_0^\infty\cdots\int_0^\infty
  \exp\Bigl[-\sum_{j=1}^nx_j^2 - \Bigl(\sum_{j=1}^nx_j\Bigr)^2\Bigr]
  dx_1\dots dx_n. 
$$
Moreover, Rogers establishes the large $n$ asymptotics
$$
  \sigma_n \sim \frac{n}{e\,2^{n/2}}.
$$
Using $\bigl(\frac{n}{n+1}\bigr)^{n/2} = \bigl((1+1/n)^n\bigr)^{-1/2}\sim
e^{-1/2}$, this yields 
$$
  \tau_n \sim \Bigl(\frac{2n}{n+1}\Bigr)^{n/2}\frac{n}{e\,2^{n/2}}
  = \Bigl(\frac{n}{n+1}\Bigr)^{n/2}\frac{n}{e}
$$
and thus the large $n$ asymptotics (cf.~\cite{Coxeter-Few-Rogers59})
$$
  \tau_n \sim \frac{n}{e^{3/2}}.
$$
Note that the exact value of $\tau_n$ is
$$
  \tau_n = \Bigl(\frac{n}{n+1}\Bigr)^{n/2}(n+1)!\,V_n
$$
In view of $\tau_n \sim \frac{n}{e^{3/2}}$ and
$\bigl(\frac{n}{n+1}\bigr)^{n/2} \sim e^{-1/2}$, this yields 
\begin{equation}\label{eq:Wn-asympt}
  W_n := n!\,V_n \sim 1/e.
\end{equation}
So to avoid rounding errors in the computation of the small quantity
$V_n$, it seems better to instead compute $W_n=n!\,V_n$ as the volume of
the region in $\R^n$ defined by the inequalities 
\begin{equation}\label{eq:Wn}
  \sum_{j=1}^nx_j^2 + \Bigl(\sum_{j=1}^nx_j\Bigr)^2\leq (n!)^{2/n}, \quad
  x_j\geq 0\text{ for }1\leq j\leq n.
\end{equation}
Using the integral expression above for $V_n$, we can also obtain
$W_n$ as
\begin{equation}\label{eq:Wn-integral}
  W_n = \frac{n!}{\Gamma(n/2+1)}\int_0^\infty\cdots\int_0^\infty
  \exp\Bigl[-\sum_{j=1}^nx_j^2 - \Bigl(\sum_{j=1}^nx_j\Bigr)^2\Bigr]
  dx_1\dots dx_n. 
\end{equation}
From this we recover $\tau_n$ as
$$
  \tau_n = \Bigl(\frac{n}{n+1}\Bigr)^{n/2}(n+1)W_n. 
$$
The following proposition, which is proved in~\S\ref{ss:Wn-bound},    
shows that the asymptotic value $1/e$ is in fact a universal lower
bound for $W_n$. 

\begin{proposition}\label{prop:Wn-bound}
The volume $W_n$ satisfies $W_n > 1/e \approx 0.368$ for all $n\in\N$. 
\end{proposition}

Using this, we obtain the following improvement of
Proposition~\ref{prop:cov-infinitesimal} in the case of ellipsoids.

\begin{proposition}\label{prop:cov-infinitesimal-ball}
Inequality~\eqref{eq:cov-infinitesimal-ball} holds for all $m\geq
n\geq 2$ except possibly for $n=2$ and $25\leq m\leq 75$. 
For $n=1$,
\eqref{eq:cov-infinitesimal-ball} holds for $m\leq 15$, $m=17$, $m=23$, $m=24$ and $m\geq 4429$.
\end{proposition}

\begin{proof}
Recall Rogers' upper bound~\eqref{eq:rogers} (specialized to balls)
$$
  \vartheta_T(B^m)\leq\wt\rho(m)\leq\rho(m),
$$
with the functions $\wt\rho,\rho$ defined above, and the simplex lower
bound $\vartheta_T(B^{m+n})\geq\tau_{m+n}\geq 1$.  
For $2\leq d\leq 24$, let $U(d)$ be the better upper bounds 
from~\cite[\S1.3]{Schuermann09}
(the covering density from Table~2 for $2\leq d\leq 5$ and from
Table~3 for $6\leq d\leq 24$, with the approximate
values rounded up in the last digit).
Set $U(1):=1$ and $U(d):=\wt\rho(d)$ for $d\geq 25$, so that
$\vartheta_T(B^d)\leq U(d)\leq\wt\rho(d)\leq\rho(d)$ for all $d\geq1$.
Then~\eqref{eq:cov-infinitesimal-ball} is implied by 
\begin{equation}\label{eq:cov-infinitesimal-binom-rogers-improved-ball}
  \frac{\tau_{m+n}}{I_{m,n,2}} > U(m)U(n),
\end{equation}
which is in turn implied by
\begin{equation}\label{eq:cov-infinitesimal-binom-rogers-ball}
  \frac{1}{I_{m,n,2}} > \rho(m)\rho(n).
\end{equation}
We will prove 
inequality~\eqref{eq:cov-infinitesimal-binom-rogers-improved-ball}
for all the pairs $m\geq n$ in the proposition.

Assume first that $n\geq 12$. We estimate $\rho(m)^2\geq\rho(m)\rho(n)$ and,
using the last assertion in Lemma~\ref{lem:I}, 
\begin{align*}
  \frac{1}{I_{m,n,2}}
  &\geq \frac{1}{I_{m,12,2}} 
  = \frac{\Gamma\Bigl(\frac{m+12}{2}+1\Bigr)}{\Gamma\Bigl(\frac{m}{2}+1\Bigr)
    \Gamma\Bigl(\frac{12}{2}+1\Bigr)} \cr
  &= \frac{(\frac{m}{2}+6)(\frac{m}{2}+5)\cdots(\frac{m}{2}+1)}{6!} 
  \geq \frac{m^6}{2^66!}.
\end{align*}
Thus~\eqref{eq:cov-infinitesimal-binom-rogers-ball} is implied by
$$
  \sigma_1(m) := \frac{\sqrt{2^66!}\rho(m)}{m^3} = \frac{96\sqrt{5}(\log m +
    \log\log m + 5)}{m^2} < 1.
$$
The function $\sigma_1$ decreases and $\sigma_1(47)\approx 0.99<1$,
so~\eqref{eq:cov-infinitesimal-binom-rogers-ball} holds if $n\geq 12$ and
$m\geq 47$.

For the remaining cases we use Proposition~\ref{prop:Wn-bound}.
From the sandwich bound
$$
  (1+1/n)^n < e < (1+1/n)^{n+1} \quad\text{for all }n\in\N
$$
we get $\bigl(\frac{n}{n+1}\bigr)^{n/2} = \bigl((1+1/n)^n\bigr)^{-1/2}
> e^{-1/2}$, so Proposition~\ref{prop:Wn-bound} implies
\begin{equation}\label{eq:taun}
  \tau_n = \Bigl(\frac{n}{n+1}\Bigr)^{n/2}(n+1)W_n >
  \frac{n+1}{e^{3/2}} \quad\text{for all }n\in\N.
\end{equation}
Assume next that $2\leq n\leq 11$.
We use the lower bound $\tau_{m+n}>(m+n+1)e^{-3/2}>m\,e^{-3/2}$
from~\eqref{eq:taun}. By the last assertion of 
Lemma~\ref{lem:I}, $1/I_{m,n,2}\geq
1/I_{m,2,2}=\frac{m}{2}+1>\frac{m}{2}$ for $n\geq 2$,
and $1/I_{m,n,2}\geq 1/I_{m,4,2}=\frac{(m+2)(m+4)}{8}>\frac{m^2}{8}$ for $n\geq 4$.
The bounds from~\cite[\S1.3]{Schuermann09} give $U(n)\leq 1.4636$ for
$n\leq 3$, and $U(n)\leq 6.3$ for $4\leq n\leq 11$.
Combining this with $U(m)\leq\rho(m)$, we conclude
that~\eqref{eq:cov-infinitesimal-binom-rogers-improved-ball} holds 
whenever $m^2/\rho(m)>2e^{3/2}\cdot 1.4636\approx 13.119$ for $n\in\{2,3\}$,
and $m^3/\rho(m)>8e^{3/2}\cdot 6.3\approx 225.88$ for $4\leq n\leq 11$.
Both $m^2/\rho(m)$ and $m^3/\rho(m)$ increase in $m$ for $m\geq 3$, and
$153^2/\rho(153)>13.13$, $49^3/\rho(49)>234$.
So~\eqref{eq:cov-infinitesimal-binom-rogers-improved-ball} 
holds for $n\in\{2,3\}$ with $m\geq 153$, and for $4\leq n\leq 11$ with $m\geq 49$.

One now checks that~\eqref{eq:cov-infinitesimal-binom-rogers-improved-ball}  
holds in the remaining finitely many cases with $m\geq n\geq 2$
($2\leq n\leq 11$ with $m\leq 152$, and $12\leq n\leq m\leq 46$)
except for the cases $n=2$ and $25\leq m\leq 75$ listed in the proposition.

Finally, let $n=1$. Note that the non-strict version
of~\eqref{eq:cov-infinitesimal-ball} holds by
Remark~\ref{rem:cov-infinitesimal-1D}. Since
$\vartheta_T(B^{m+1})\geq\tau_{m+1}$, the strict inequality holds
whenever $U(m)<\tau_{m+1}/I_{m,1,2}$. For $m=1$ this follows from
$\tau_2\geq 1$ and $1/I_{1,1,2}=4/\pi>1=U(1)$. For $2\leq m\leq 15$,
$m=17$, $m=23$ and $m=24$ it follows by direct computation from
$\tau_{m+1}\geq\max\{1,(m+2)e^{-3/2}\}$ and the values $U(m)$
in~\cite[\S1.3]{Schuermann09}.
For $4429\leq m\leq 8305$ it follows by direct computation from
$\tau_{m+1}>(m+2)e^{-3/2}$ and $U(m)=\wt\rho(m)$. For $m\geq 8306$,
log-convexity of $\Gamma$ gives $1/I_{m,1,2}\geq\sqrt{2(m+1)/\pi}$,
hence $\tau_{m+1}/I_{m,1,2}>c\,m^{3/2}$ with $c=\sqrt{2/\pi}\,e^{-3/2}$,
while $U(m)\leq\rho(m)$. Since $m^{3/2}/\rho(m)$ increases in $m$ and
$c\,8306^{3/2}>\rho(8306)$, the strict inequality follows.
\end{proof}

\begin{remark}
The remaining cases with $n=2$ are off by a factor of at most $2.2$, so to cover
them one would need to improve some bound by a factor of $2.2$. 
\end{remark}

\begin{remark}
Numerical computations
suggest that the volume $W_n$ strictly decreases in $n\in\N$ from $W_1\approx
0.707$ to $\lim_{n\to\infty}W_n = 1/e \approx 0.368$. 
\end{remark}

\section{Finsler manifolds}\label{sec:Finsler}

A (continuous)\footnote{
Often the definition of a Finsler metric includes smoothness
conditions, which will not be needed in this note.}
{\em Finsler metric} on a manifold $M$ is a continuous
function $F:TM\to\R_{\geq 0}$ satisfying for all $x\in M$ and $v,w\in T_xM$:
\begin{itemize}
\item $F_x(v+w)\leq F_x(v)+F_x(w)$;
\item $F_x(\lambda v)=\lambda F_x(v)$ for all $\lambda\geq 0$;
\item $F_x(v)>0$ for $v\neq 0$.
\end{itemize}
In other words, $F$ gives a continuous family of (possibly
nonsymmetric) norms $F_x$ on the tangent spaces.
By a {\em path from $x$ to $y$} we will mean a $C^1$-map
$\gamma:[0,1]\to M$ with $\gamma(0)=x$ and $\gamma(1)=y$. Its length
with respect to $F$ is
$$
  L_F(\gamma) := \int_0^1F_{\gamma(t)}\bigl(\dot\gamma(t)\bigr)dt. 
$$
The corresponding (possibly nonsymmetric) distance from $x$ to $y$ is
$$
  d_F(x,y) := \inf\{L_F(\gamma)\mid \gamma\text{ path from $x$ to $y$}\}.
$$
When $F$ is clear from the context we will omit it and denote $d_F$ by
$d_M$. 

For two Finsler manifolds $M,N$ and $p\in[1,\infty]$, the $p$-norm on
the tangent space defines a Finsler metric on $M\times N$.
We define the (possibly nonsymmetric) metric $d_p$ on $M\times N$ 
by formula~\eqref{eq:dp} with the metrics $d_M,d_N$.  

The following two theorems correspond to
Theorems~\ref{thm:Finsler-intro} and~\ref{thm:Riemann-intro} in the
Introduction. 

\begin{theorem}\label{thm:Finsler}
For each pair of compact Finsler manifolds $M,N$ with $m=\dim M,
n=\dim N$ as in Proposition~\ref{prop:cov-infinitesimal} there exist
constants $r_{M,N}>0$ and $p_{m,n}>1$ such that for all $0<r\leq r_{M,N}$
and $1\leq p\leq p_{m,n}$ we have
$$
  \cov(M\times N,d_p,r) > \cov(M,d_M,r)\,\cov(N,d_N,r).
$$
\end{theorem}

\begin{theorem}\label{thm:Riemann}
For each pair of compact Riemannian manifolds $M,N$ with $m=\dim M,
n=\dim N$ as in Proposition~\ref{prop:cov-infinitesimal-ball}
there exists a constant $r_{M,N}>0$ such that for all $0<r\leq r_{M,N}$ we have
$$
  \cov(M\times N,d_2,r) > \cov(M,d_M,r)\,\cov(N,d_N,r).
$$
\end{theorem}

The remainder of this section is devoted to the proof of these
theorems. We begin with some preparation. 

{\bf Estimates on Finsler metrics. }
Consider a Finsler metric $F$ on $\R^m$ and its nonsymmetric norms
$F_x$, $x\in\R^m$. By continuity of $F$, for each $\lambda>1$ there
exists $\rho>0$ such that 
\begin{equation}\label{eq:norm-comparison}
  \lambda^{-1}F_0(v)\leq F_x(v)\leq \lambda F_0(v)\quad\text{for all
    $x,v\in\R^m$ with $F_0(x)\leq 5\rho$}. 
\end{equation}
The first lemma compares the balls
$$
  B_{F_0}(x,r) = \{y\in\R^m\mid F_0(y-x)\leq r\},\quad 
  B_{d_F}(x,r) = \{y\in\R^m\mid d_F(x,y)\leq r\}.
$$

\begin{lemma}[Constant approximation]\label{lem:constant-approx}
For each $\lambda>1$ there exists $\rho>0$ such that for all
$x\in\R^m$ with $F_0(x)\leq \rho$ and $0<r\leq\rho$ we have
\begin{equation}\label{eq:constant-approx}
  B_{F_0}(x,\lambda^{-1}r) \subset B_{d_F}(x,r) \subset B_{F_0}(x,\lambda r).
\end{equation}
\end{lemma}

\begin{proof}
Both inclusions in~\eqref{eq:constant-approx} become weaker
as $\lambda$ increases, so (replacing $\lambda$ by $\min\{\lambda,2\}$)
we may assume $\lambda\leq 2$.
Given $\lambda>1$, we take $\rho>0$ such
that~\eqref{eq:norm-comparison} holds. 
For the first inclusion, we need to show
$$
  d_F(x,y) \leq \lambda F_0(y-x)\text{ for all $x,y\in\R^m$ with
    $F_0(x)\leq\rho$, $F_0(y-x)\leq\lambda^{-1}\rho$}.
$$
To see this, set $v:=y-x$ and take the path $\gamma(t):=x+tv$. Then
using~\eqref{eq:norm-comparison} we get
$$
  d_F(x,y) \leq L_F(\gamma)
  = \int_0^1 F_{\gamma(t)}(v)dt
  \leq \int_0^1 \lambda F_0(v)dt
  = \lambda F_0(v). 
$$
For the second inclusion, we will need the following fact: for each $C^1$-path $\gamma:[a,b]\to\R^m$,
\begin{equation}\label{eq:norm-integral}
   F_0\bigl(\gamma(b)-\gamma(a)\bigr) \leq \int_a^b F_0\bigl(\dot\gamma(t)\bigr)dt.
\end{equation}
To see this, take for large $N$ the subdivision $t_i=a+i(b-a)/N$ and write
$$
  \gamma(b)-\gamma(a) = \sum_{i=1}^N\frac{\gamma(t_i)-\gamma(t_{i-1})}{t_i-t_{i-1}}(t_i-t_{i-1}).
$$
Applying $F_0$ and using the triangle inequality, we obtain
\begin{align*}
  F_0\bigl(\gamma(b)-\gamma(a)\bigr)
  &\leq \sum_{i=1}^NF_0\Bigl(\frac{\gamma(t_i) -
      \gamma(t_{i-1})}{t_i-t_{i-1}}\Bigr)(t_i-t_{i-1}) \cr 
  &\approx \sum_{i=1}^NF_0\bigl(\dot\gamma(t_i)\bigr)(t_i-t_{i-1})  
  \approx \int_a^bF_0\bigl(\dot\gamma(t)\bigr)dt, 
\end{align*}
where $\approx$ means equality up to an error that vanishes as $N\to\infty$.

Suppose now that $F_0(x)\leq\rho$, $0<r\leq\rho$, and $y\in B_{d_F}(x,r)$.
Then for each $\eps>0$ there exists a path $\gamma$ from $x$ to $y$
with $L_F(\gamma)\leq r+\eps$. For $b\in[0,1]$,
using~\eqref{eq:norm-integral} and~\eqref{eq:norm-comparison} we compute
$$
  F_0\bigl(\gamma(b)-x\bigr)
  \leq \int_0^b F_0\bigl(\dot\gamma(t)\bigr)dt
  \leq \lambda\int_0^b F_{\gamma(t)}\bigl(\dot\gamma(t)\bigr)dt
  \leq \lambda L_F(\gamma)
  \leq \lambda(r+\eps).
$$
We may assume without loss of generality that $\eps<\rho$. 
Let $b^*$ be the supremum of all $b\in[0,1]$ such that
$F_0\bigl(\gamma(t)\bigr)\leq 5\rho$ for all $t\in[0,b]$; by continuity
this also holds for $b=b^*$. For $b\leq b^*$ the preceding
computation is justified and gives
$F_0\bigl(\gamma(b)\bigr)\leq F_0(x)+F_0\bigl(\gamma(b)-x\bigr)\leq
\rho+\lambda(r+\eps)<\rho+4\rho=5\rho$ (where we have used
$\lambda\leq 2$). If $b^*<1$, then continuity would
give $F_0\bigl(\gamma(t)\bigr)\leq 5\rho$ slightly beyond $b^*$, a
contradiction. So $b^*=1$, which justifies the use
of~\eqref{eq:norm-comparison} in the preceding computation for all $b\in[0,1]$.
For $b=1$ and $\eps\to 0$, the last displayed formula becomes
$F_0(y-x)\leq\lambda r$ and the second inclusion of the lemma follows.
\end{proof}
 
{\bf Estimates on covering numbers. }
Consider a convex body $A\subset\R^m$ and a compact subset
$X\subset\R^m$ satisfying~\eqref{eq:X}. 

\begin{lemma}[Rescaling]\label{lem:rescaling}
(a) For all $\lambda\geq 1$ and $r>0$ we have
$$
  \cov(X,d_A,\lambda r) \leq \cov(X,d_A,r). 
$$
(b) For each $\lambda>1$ there exists $r_{X,A,\lambda}>0$ such that
  for all $0<r\leq r_{X,A,\lambda}$ we have
$$
  \cov(X,d_A,r) \leq \lambda^{m+2}\cov(X,d_A,\lambda r).
$$
\end{lemma}

\begin{proof}
Part (a) is clear. For (b), fix $\lambda>1$. By
Proposition~\ref{prop:Groemer} and Lemma~\ref{lem:cov-ext}, there
exists $r_{X,A,\lambda}>0$ 
such that for all $0<s\leq\lambda r_{X,A,\lambda}$ we have
\begin{gather*}
  \lambda^{-1/2}\vartheta_T(A) \leq \vartheta_T(X;sA) \leq \lambda^{1/2}\vartheta_T(A)
  \quad\text{and}\cr
  \cov(X,d_A,s)\leq\lambda\,\cov^{\rm ext}(X,d_A,s).
\end{gather*}
For $0<r\leq r_{X,A,\lambda}$ we have $r,\lambda r\leq\lambda
r_{X,A,\lambda}$, and using~\eqref{eq:cov-number-density} we estimate
\begin{align*}
  \cov(X,d_A,r)
  &\leq \lambda\cov^{\rm ext}(X,d_A,r)
  = \lambda\frac{\vol(X)}{\vol(rA)}\vartheta_T(X;rA) \cr
  &\leq \lambda^{3/2}\frac{\vol(X)}{\vol(rA)}\vartheta_T(A) \cr
  &= \lambda^{m+3/2}\frac{\vol(X)}{\vol(\lambda rA)}\vartheta_T(A)
  \leq \lambda^{m+2}\frac{\vol(X)}{\vol(\lambda rA)}\vartheta_T(X;\lambda rA) \cr
  &= \lambda^{m+2}\cov^{\rm ext}(X,d_A,\lambda r)
  \leq \lambda^{m+2}\cov(X,d_A,\lambda r).
\end{align*}
\end{proof}

\begin{lemma}[Volume bounds]\label{lem:vol-bounds}
Let $F$ be a Finsler metric on $\R^m$ and $X\subset\R^m$ a compact
subset satisfying~\eqref{eq:X}. 
There exist constants $r_X>0$ and $C_F\geq 1$ such that for all
$0<r\leq r_X$ we have  
$$
  \frac{1}{C_F}\frac{\vol(X)}{r^m} \leq \cov(X,d_F,r) \leq
  C_F\frac{\vol(X)}{r^m}. 
$$
The same constant $C_F$ works for every compact subset $Y\subset X$
satisfying~\eqref{eq:X} in place of $X$, with some threshold $r_Y>0$
in place of $r_X$. 
\end{lemma}

\begin{proof}
Let $Q=[-1/2,1/2]^m$ be the unit cube centered at $0$. 
By compactness of $X$, there exist constants $K\geq 1$
and $r'>0$ depending on
$F$ and $X$ such that for all $x\in X$ and $0<r\leq r'$ we have
$$
  x+\frac{r}{K}Q \subset B_{d_F}(x,r) \subset x+KrQ.
$$
Since $\vartheta_T(Q)=1$, by Proposition~\ref{prop:Groemer} and
Lemma~\ref{lem:cov-ext} there
exists $r_X\in(0,r']$ such that for all $0<r\leq r_X$ we have
$1/2\leq \vartheta_T(X;rQ)\leq 2$ and $\cov(X,d_Q,r)\leq
2\cov^{\rm ext}(X,d_Q,r)$, and therefore
\begin{align*}
  \cov(X,d_F,r) &\leq \cov(X,d_Q,r/K)
  \leq 2\cov^{\rm ext}(X,d_Q,r/K) \cr
  &= \frac{2\vol(X)}{(r/K)^m}\vartheta_T(X;\frac{r}{K}Q) \cr
  &\leq \frac{4\vol(X)}{(r/K)^m} 
  = 4K^m\frac{\vol(X)}{r^m}. 
\end{align*}
This proves the second estimate with $C_F=4K^m$, and the first one
follows analogously. 
Since $K$, $r'$ and $C_F$ depend only on $F$ and $X$, and the cube
comparison holds at all points of $Y\subset X$, the same argument applies to
every compact $Y\subset X$ satisfying~\eqref{eq:X}, with $r_Y\in(0,r']$ chosen
for $Y$.
\end{proof}

\begin{lemma}[Shrinking]\label{lem:shrinking}
Let $F$ be a Finsler metric on $\R^m$ and $X\subset\R^m$ a convex
polyhedron with $0\in\Int X$. Then for each $\mu>1$ there exist
$\lambda\in(0,1)$ and $\rho>0$ 
(depending on $X,\mu$) such that for all $0<r\leq\rho$ we have
$$
   \cov(\lambda X;X,d_F,r)\geq \mu^{-1}\cov(X,d_F,r).
$$
Here $\cov(\lambda X;X,d_F,r)$ counts balls with centers in $X$
covering $\lambda X$ as in~\S\ref{sec:basic}.
\end{lemma}

\begin{proof}
For $\lambda\in(0,1)$ let $X'=\lambda X$. Then
$X'\subset X$ since $X$ is convex and $0\in X$. 
By Lemma~\ref{lem:vol-bounds}, applied to $X$ and to the compact
subset $\overline{X\setminus X'}\subset X$ (which
satisfies~\eqref{eq:X}), there exist constants $r_{X,\lambda}>0$ 
and $C_F\geq 1$ such that for all $0<r\leq r_{X,\lambda}$ we have  
$$
  \frac{1}{C_F}\frac{\vol(X)}{r^m} \leq \cov(X,d_F,r)
  \quad\text{and}\quad
  \cov(\overline{X\setminus X'},d_F,r) \leq C_F\frac{\vol(X\setminus X')}{r^m}. 
$$
(Note that $\vol(\overline{X\setminus X'})=\vol(X\setminus X')$ because
$\p X'$ is a null set.) Since every point of $X$ lies in $X'$ or in
$\overline{X\setminus X'}\subset X$, we have $\cov(X,d_F,r)\leq
\cov(X';X,d_F,r)+\cov(\overline{X\setminus X'},d_F,r)$.
Using this, we estimate
\begin{align*}
  1 - \frac{\cov(X';X,d_F,r)}{\cov(X,d_F,r)}
  &\leq \frac{\cov(\overline{X\setminus X'},d_F,r)}{\cov(X,d_F,r)}  
  \leq \frac{C_F\vol(X\setminus X')/r^m}{\vol(X)/C_Fr^m} \cr 
  &= C_F^2\frac{\vol(X\setminus X')}{\vol(X)}. 
\end{align*}
Since the constant $C_F$ is the one of Lemma~\ref{lem:vol-bounds}
for $X$, which also works for all compact subsets of $X$, it does not
depend on $\lambda$, and we  can choose $\lambda$ such that
$$
  \frac{\cov(X';X,d_F,r)}{\cov(X,d_F,r)}
  \geq 1 - C_F^2\frac{\vol(X\setminus X')}{\vol(X)} \geq \mu^{-1}.
$$
\end{proof}

{\bf Localization. }
The following lemma reduces the setup of Theorem~\ref{thm:Finsler} to
simplices in Euclidean space.  

\begin{lemma}[Localization]\label{lem:localization}
Consider compact Finsler manifolds $M,N$ with triangulations
$M=\bigcup_iM_i$ and $N=\bigcup_jN_j$. Let $p\geq 1$. Suppose that
there exist $\mu>1$ and $r_0>0$ such that 
$$
  \cov(M_i\times N_j,d_p,r) \geq \mu^2\cov(M_i,d_M,r)\,\cov(N_j,d_N,r)
$$
for all $i,j$ and $0<r\leq r_0$. Then there exists $r_1\in(0,r_0)$ such
that  
$$
  \cov(M\times N,d_p,r) \geq \mu\,\cov(M,d_M,r)\,\cov(N,d_N,r)
$$
for all $0<r\leq r_1$. 
\end{lemma}

\begin{proof}
For each $i$, we pick a diffeomorphism of an open neighbourhood of
$M_i$ onto $\R^m$ sending $M_i$ to a simplex containing the origin in
its interior, and we pull back the Finsler metric from $M$ to a Finsler
metric on $\R^m$.
Similarly, we embed each $N_j$ into $\R^n$ and pull back the Finsler metric. 
We identify $M_i$ and $N_j$ with their images in $\R^m$ and $\R^n$
and apply Lemma~\ref{lem:shrinking} to the convex polytopes $M_i\times
N_j$ with the metric $d_p$. This yields $\lambda<1$ and $\rho>0$ such
that for all $i,j$ and all $0<r\leq\rho$ we have 
$$
   \cov\bigl(\lambda(M_i\times N_j);M_i\times N_j,d_p,r\bigr)\geq
   \mu^{-1}\cov(M_i\times N_j,d_p,r). 
$$
Pick $0<r_1\leq\min\{r_0,\rho\}$ so small that for all
$i,j$ each $r_1$-ball with respect to the metric $d_p$
intersecting $\lambda (M_i\times N_j)$ is contained in $\Int(M_i\times N_j)$.
Then for $0<r\leq r_1$ each covering of $M\times N$ by $r$-balls
contains non-overlapping coverings of $\bigcup_{i,j}\lambda(M_i\times N_j)$
(with centers in $M_i\times N_j$)
and we obtain  
\begin{align*}
  \cov(M\times N,d_p,r)
  &\geq \sum_{i,j}\cov\bigl(\lambda(M_i\times N_j);M_i\times N_j,d_p,r\bigr) \cr
  &\geq \mu^{-1}\sum_{i,j}\cov(M_i\times N_j,d_p,r) \cr
  &\geq \mu\sum_{i,j}\cov(M_i,d_M,r)\,\cov(N_j,d_N,r) \cr
  &\geq \mu\,\cov(M,d_M,r)\,\cov(N,d_N,r).
\end{align*}
\end{proof}

{\bf Proof of the theorems. }
Now we are ready to prove the two theorems.

\begin{proof}[Proof of Theorem~\ref{thm:Finsler}]
Consider compact Finsler manifolds $M,N$ with $m=\dim M,
n=\dim N$ as in Proposition~\ref{prop:cov-infinitesimal}. 
By Corollary~\ref{cor:cov-infinitesimal}, there exist constants
$\mu=\mu_{m,n}>1$ and $p_{m,n}>1$ depending only on
$m,n$ such that for all convex bodies $A\subset\R^m$, $B\subset\R^n$
and $1\leq p\leq p_{m,n}$ we have
\begin{equation*}
  \frac{\vartheta_T(A\times_pB)}{I_{m,n,p}} \geq \mu^4\vartheta_T(A)\vartheta_T(B).
\end{equation*}
Since $(a^p+b^p)^{1/p}$ is nonincreasing in $p$, we have $d_p\geq d_{p'}$
for $1\leq p\leq p'$, so every $d_p$-ball is contained in the
$d_{p'}$-ball with the same center and radius, and
$\cov(M\times N,d_p,r)\geq\cov(M\times N,d_{p_{m,n}},r)$ for all
$p\in[1,p_{m,n}]$. It therefore suffices to prove the theorem for
$p=p_{m,n}$, which we fix from now on.
Fix a constant $\lambda>1$ satisfying
\begin{equation}\label{eq:lambda}
  \lambda^{2m+2n+6}\leq\mu. 
\end{equation}
Denote the Finsler metric on $M$ by $F$, so that $d_M=d_F$.
Consider a triangulation $M=\bigcup_iM_i$. As in the proof of
Lemma~\ref{lem:localization}, we identify each $M_i$ with a simplex in 
$\R^m$ containing the origin in its interior. 
By Lemma~\ref{lem:constant-approx} (after translating $z$ to the origin),
for each $z\in M_i$ there exists $\rho_z>0$ such that for all
$x\in M_i$ with $F_z(x-z)\leq \rho_z$ and $0<r\leq\rho_z$ we have
\begin{equation*}
  B_{F_z}(x,\lambda^{-1}r) \subset B_{d_F}(x,r) \subset B_{F_z}(x,\lambda r).
\end{equation*}
The sets 
$$
  U_z := \{x\in M_i\mid F_z(x-z)< \rho_z\},\qquad z\in M_i
$$
form an open covering of $M_i$. By compactness, after subdividing the
$M_i$, we may assume that each $M_i$ is contained in one such
set. That is, there exist $z_i\in M_i$ and $\rho_i>0$ such that for
all $x\in M_i$ and $0<r\leq\rho_i$ we have
\begin{equation*}
  B_{F_{z_i}}(x,\lambda^{-1}r) \subset B_{d_F}(x,r) \subset B_{F_{z_i}}(x,\lambda r).
\end{equation*}
With the convex bodies $A_i=B_{F_{z_i}}(0,1)\subset\R^m$ and
$d_F=d_M$,
this reads
\begin{equation*}
  x+\lambda^{-1}rA_i \subset B_{d_M}(x,r) \subset x+\lambda rA_i.
\end{equation*}
Similarly, after subdividing the $N_j$, there exist $\sigma_j>0$ and
convex bodies $B_j\subset\R^n$ such that for all $y\in N_j$ and
$0<r\leq\sigma_j$ we have 
\begin{equation*}
  y+\lambda^{-1}rB_j \subset B_{d_N}(y,r) \subset y+\lambda rB_j.
\end{equation*}
For $0<r\leq\rho:=\min\{\rho_i,\sigma_j\}$ this implies the following
inequalities for the covering numbers:
\begin{gather*}
  \cov(M_i,A_i,\lambda r) \leq \cov(M_i,d_M,r) \leq \cov(M_i,A_i,\lambda^{-1}r),\cr
  \cov(N_j,B_j,\lambda r) \leq \cov(N_j,d_N,r) \leq \cov(N_j,B_j,\lambda^{-1}r),\cr
  \cov(M_i\times N_j,A_i\hspace{-2pt}\times_pB_j,\lambda r) \leq \cov(M_i\times
  N_j,d_p,r) \leq \cov(M_i\times N_j,A_i\hspace{-2pt}\times_pB_j,\lambda^{-1}r).
\end{gather*}
By Corollary~\ref{cor:cov-infinitesimal}, 
there exists $r_0>0$ such that for all $i,j$ and 
$0<r\leq r_0$ we have
\begin{equation*}
  \cov(M_i\times N_j,A_i\hspace{-2pt}\times_pB_j,r) \geq \mu^3\cov(M_i,A_i,r)\,\cov(N_j,B_j,r).
\end{equation*}
Using the preceding inequalities together with
Lemma~\ref{lem:rescaling} (for the fixed $\lambda$, with $r$ below
its thresholds $r_{X,A,\lambda}$ for the finitely many sets $M_i$, $N_j$,
$M_i\times N_j$) and~\eqref{eq:lambda}, we estimate
\begin{align*}
  \frac{\cov(M_i\times N_j,d_p,r)}{\cov(M_i,d_M,r)\cov(N_j,d_N,r)} 
  &\geq \frac{\cov(M_i\times N_j,A_i\hspace{-2pt}\times_pB_j,\lambda
    r)}{\cov(M_i,A_i,\lambda^{-1}r)\cov(N_j,B_j,\lambda^{-1}r)} \cr
  &\geq \frac{\lambda^{-(m+n+2)}\cov(M_i\times N_j,A_i\hspace{-2pt}\times_pB_j,
    r)}{\lambda^{m+2}\cov(M_i,A_i,r)\lambda^{n+2}\cov(N_j,B_j,r)} \cr
  &\geq \frac{\mu^3}{\lambda^{2m+2n+6}}
  \geq \mu^2. 
\end{align*}
The theorem now follows from Lemma~\ref{lem:localization}.
\end{proof}

\begin{proof}[Proof of Theorem~\ref{thm:Riemann}]
This proof is entirely analogous to the proof of
Theorem~\ref{thm:Finsler}, using
Proposition~\ref{prop:cov-infinitesimal-ball} instead of 
Proposition~\ref{prop:cov-infinitesimal}. 
\end{proof}

\appendix

\section{Some estimates}\label{sec:estimates}

In this appendix we collect some estimates that are used in the paper.

\subsection{Sandwich bounds}

\begin{lemma}[Sandwich bounds for the Euler number]\label{lem:sandwich-euler}
For each real number $x\geq 1$ we have
$$
  \Bigl(1+\frac{1}{x}\Bigr)^x < e < \Bigl(1+\frac{1}{x}\Bigr)^{x+1}
$$
\end{lemma}

\begin{proof}
Recall that the derivative of the function $a^x$ is $a^x\log a$, and
the power series expansion of the logarithm yields $y-y^2/2<\log(1+y)<y$
for all $y\in(0,1)$, which we will apply with $y=1/x$. 
Consider the functions
$$
  f(x) := \Bigl(1+\frac{1}{x}\Bigr)^x,\qquad
  g(x) := \Bigl(1+\frac{1}{x}\Bigr)^{x+1}
$$
and compute for $x>1$:
\begin{align*}
  f'(x)
  &= \Bigl(1+\frac{1}{x}\Bigr)^x\log\Bigl(1+\frac{1}{x}\Bigr)
  + x\Bigl(1+\frac{1}{x}\Bigr)^{x-1}\Bigl(-\frac{1}{x^2}\Bigr) 
  = \Bigl(1+\frac{1}{x}\Bigr)^{x-1}\phi(x) \cr
  \phi(x)
  &= \Bigl(1+\frac{1}{x}\Bigr)\log\Bigl(1+\frac{1}{x}\Bigr)-\frac{1}{x} 
  > \Bigl(1+\frac{1}{x}\Bigr)\Bigl(\frac{1}{x}-\frac{1}{2x^2}\Bigr)-\frac{1}{x} \cr
  &= \frac{1}{2x^2}-\frac{1}{2x^3} > 0, \cr
  g'(x)
  &= \Bigl(1+\frac{1}{x}\Bigr)^{x+1}\log\Bigl(1+\frac{1}{x}\Bigr)
  + (x+1)\Bigl(1+\frac{1}{x}\Bigr)^x\Bigl(-\frac{1}{x^2}\Bigr) 
  = \Bigl(1+\frac{1}{x}\Bigr)^x\psi(x) \cr
  \psi(x)
  &= \Bigl(1+\frac{1}{x}\Bigr)\log\Bigl(1+\frac{1}{x}\Bigr)-\frac{x+1}{x^2} 
  < \Bigl(1+\frac{1}{x}\Bigr)\frac{1}{x} - \frac{x+1}{x^2} = 0.
\end{align*}
Now the lemma follows from
$\lim_{x\to\infty}f(x)=\lim_{x\to\infty}g(x)=e$, since $f$ is strictly
increasing and $g$ strictly decreasing. 
\end{proof}

\begin{lemma}[Sandwich bounds for the Gamma function]\label{lem:sandwich-gamma}
For each real number $x\geq 1$ we have
$$
  \sqrt{\frac{x}{x+1}}\sqrt{2\pi x}\Bigl(\frac{x}{e}\Bigr)^x <
  \Gamma(x+1) <
  \sqrt{\frac{x+1}{x}}e^{\frac{1}{12(x+1)}}\sqrt{2\pi x}\Bigl(\frac{x}{e}\Bigr)^x.  
$$
\end{lemma}

\begin{proof}
For each $z>0$ we have the well-known bounds (see e.g.~\cite{Nemes15})
$$
  1 < \frac{\Gamma(z)}{\sqrt{2\pi}z^{z-1/2}e^{-z}} < e^{\frac{1}{12z}}.
$$
For $z=x+1$ with $x\geq 1$ we use Lemma~\ref{lem:sandwich-euler} to
estimate the denominator from below and above:   
\begin{align*}
  \sqrt{2\pi}(x+1)^{x+1/2}e^{-x-1}
  &= \sqrt{\frac{2\pi}{x+1}}\Bigl(1+\frac{1}{x}\Bigr)^{x+1}
  \Bigl(\frac{x}{e}\Bigr)^{x+1}\cr 
  &> \sqrt{\frac{2\pi}{x+1}}e \Bigl(\frac{x}{e}\Bigr)^{x+1} 
  = \sqrt{\frac{x}{x+1}} \sqrt{2\pi x}\Bigl(\frac{x}{e}\Bigr)^x,\cr 
  \sqrt{2\pi}(x+1)^{x+1/2}e^{-x-1}
  &= \frac{\sqrt{2\pi(x+1)}}{e}\Bigl(1+\frac{1}{x}\Bigr)^x
  \Bigl(\frac{x}{e}\Bigr)^x\cr 
  &< \frac{\sqrt{2\pi(x+1)}}{e}e\Bigl(\frac{x}{e}\Bigr)^x
  = \sqrt{\frac{x+1}{x}} \sqrt{2\pi x}\Bigl(\frac{x}{e}\Bigr)^x.
\end{align*}
\end{proof}

\subsection{A lower estimate on Rogers' simplex bound}\label{ss:Wn-bound}

Here we prove Proposition~\ref{prop:Wn-bound} providing a lower bound
on the volume $W_n$ in~\ref{eq:Wn-asympt}, 
and thus on Rogers' simplex bound $\tau_n$. 
This subsection was produced with the help of GPT-5.2-pro. 

For $n\in\mathbb{N}$, recall that $W_n$ is the $n$-dimensional volume
of the region 
\begin{align*}
  \Omega_n &:= \Bigl\{x=(x_1,\dots,x_n)\in\mathbb{R}^n\;\Bigl|\; x_i\ge
  0,\ \sum_{i=1}^n x_i^2+\Bigl(\sum_{i=1}^n x_i\Bigr)^2\le
  (n!)^{2/n}\Bigr\} \cr
  &= \{x\in\R^n \mid x_i\geq 0,\ Q(x)+S(x)^2\le R_n^2\}
\end{align*}
with
\[
  R_n:=(n!)^{1/n},\qquad S(x):=\sum_{i=1}^n x_i,\qquad Q(x):=\sum_{i=1}^n x_i^2.
\]
We will first write $W_n$ as an average over the standard simplex
\[
\Delta_{n-1}:=\Bigl\{u\in\mathbb{R}^n\;\Bigl|\; u_i\ge 0,\ \sum_{i=1}^n u_i=1\Bigr\}
\]
with its usual $(n-1)$-dimensional Lebesgue measure $du=du_2\cdots du_n$.


\begin{lemma}[Simplicial radial coordinates]\label{lem:jac}
Every $x\in\mathbb{R}^n_{\ge 0}\setminus\{0\}$ can be written uniquely as $x=s u$ with
\[
s=S(x)>0,\qquad u=\frac{x}{S(x)}\in\Delta_{n-1}.
\]
In coordinates $(s,u_2,\dots,u_n)$ the Jacobian is $s^{n-1}$, i.e.
\[
dx = s^{n-1}\,ds\,du.
\]
\end{lemma}

\begin{proof}
Uniqueness is clear from $s=S(x)$ and $u=x/s$. The Jacobian computation is standard:
write $x_i=s u_i$ for $2\le i\le n$ and $x_1=s(1-\sum_{i=2}^{n}u_i)$ and compute the determinant;
it equals $s^{n-1}$.
\end{proof}

If $x=s u$ with $u\in\Delta_{n-1}$, then
\[
S(x)=s,\qquad Q(x)=s^2\sum_{i=1}^n u_i^2,
\]
hence the constraint $Q(x)+S(x)^2\le R_n^2$ becomes
\[
s^2\Bigl(1+\sum_{i=1}^n u_i^2\Bigr)\le R_n^2
\quad\Longleftrightarrow\quad
0\le s\le \frac{R_n}{\sqrt{1+\sum_{i=1}^n u_i^2}}.
\]
Therefore, by Lemma~\ref{lem:jac},
\begin{align}
W_n
&=\int_{\Delta_{n-1}}\int_{0}^{R_n/\sqrt{1+\sum u_i^2}} s^{n-1}\,ds\,du\nonumber\\
&=\int_{\Delta_{n-1}}\frac1n\left(\frac{R_n}{\sqrt{1+\sum u_i^2}}\right)^n\,du\nonumber\\
&=\frac{R_n^n}{n}\int_{\Delta_{n-1}}\left(1+\sum_{i=1}^n u_i^2\right)^{-n/2}\,du.\label{eq:Wn-simplex}
\end{align}
Since $R_n^n=n!$, we have $\frac{R_n^n}{n}=(n-1)!$.
Also $\vol(\Delta_{n-1})=\frac1{(n-1)!}$. Hence, if we define the \emph{simplex average}
\[
\langle f\rangle_{\Delta}:=\frac{1}{\vol(\Delta_{n-1})}\int_{\Delta_{n-1}} f(u)\,du=(n-1)!\int_{\Delta_{n-1}} f(u)\,du,
\]
then \eqref{eq:Wn-simplex} becomes the identity
\begin{equation}\label{eq:Wn-average}
W_n=\left\langle\left(1+\sum_{i=1}^n u_i^2\right)^{-n/2}\right\rangle_{\Delta}.
\end{equation}


%

Define for $u\in\Delta_{n-1}$ the quantity
\[
  T(u):=\sum_{i=1}^n u_i^2.
\]

\begin{lemma}
\label{lem:Tmoments}
The simplex average of $T$ equals
\[
\langle T\rangle_{\Delta}=\frac{2}{n+1}.
\]
\end{lemma}

\begin{proof}
By symmetry, $\langle T\rangle_{\Delta}=n\langle u_n^2\rangle_{\Delta}$.
Using the beta function identity $\int_0^1
t^p(1-t)^q\,dt=\frac{p!\,q!}{(p+q+1)!}$ for integers $p,q\ge 0$ (which
is easily proved by partial integration and induction), we compute
\begin{align*}
  \frac{\langle T\rangle_{\Delta}}{n!}
  &= \frac{\langle u_n^2\rangle_{\Delta}}{(n-1)!} 
  = \int_{\Delta_{n-1}} u_n^2du_2\cdots du_n 
  = \int_0^1 u_n^2\frac{(1-u_n)^{n-2}}{(n-2)!}du_n \cr
  &= \frac{1}{(n-2)!}\frac{2!\,(n-2)!}{(n+1)!} 
  = \frac{2}{(n+1)!}.
\end{align*}
\end{proof}

Now we are ready to prove Proposition~\ref{prop:Wn-bound} in the
following more precise form.

\begin{proposition}\label{prop:Wn-bound-precise}
For every $n\in\mathbb{N}$ we have
\[
W_n \ge \left(1+\frac{2}{n+1}\right)^{-n/2} > \frac1e.
\]
\end{proposition}

\begin{proof}
Let $\phi(t):=(1+t)^{-n/2}$ for $t\ge 0$. Then
\[
\phi''(t)=\frac{n}{2}\left(\frac{n}{2}+1\right)(1+t)^{-n/2-2}>0,
\]
so $\phi$ is convex. Jensen's inequality in integral form (see
e.g.~\cite{Rudin87}) gives
\[
\left\langle \phi(T)\right\rangle_{\Delta} \ge \phi\!\left(\langle T\rangle_{\Delta}\right).
\]
By \eqref{eq:Wn-average} and Lemma~\ref{lem:Tmoments},
\[
W_n=\langle \phi(T)\rangle_{\Delta}\ge \phi\!\left(\frac{2}{n+1}\right)
=\left(1+\frac{2}{n+1}\right)^{-n/2}.
\]
To compare with $e^{-1}$,
apply Lemma~\ref{lem:sandwich-euler} with $x=(n+1)/2\geq 1$ to get
\[
\left(1+\frac{2}{n+1}\right)^{(n+1)/2} < e.
\]
Raising both sides to the power $-\frac{n}{n+1}$ yields
\[
\left(1+\frac{2}{n+1}\right)^{-n/2} > e^{-n/(n+1)} > e^{-1}.
\]
Thus $W_n>e^{-1}$ for every $n$.
\end{proof}

\section{Covering numbers of random metric spaces}\label{sec:random}

Together with M.\,Zerbin we have probed Conjecture~\ref{conj:cov} on
random metric spaces (without finding a counterexample). In this
appendix we briefly describe the underlying algorithms. 

{\bf Random metric spaces. }
Here we describe Vershik's algorithm~\cite[\S2.4]{Vershik04} for
generating random metrics on the set $\{1,\dots,n\}$ for $n\geq 2$. 
The metric is described by the matrix $R$ of distances
$r_{i,j}=d(i,j)$, $1\leq i,j\leq n$, satisfying
$$
  r_{i,j}=r_{j,i}\geq 0,\qquad r_{i,i}=0
$$
and the triangle inequalities. Note that it suffices to specify the
upper triangular part $r_{i,j}$ for $i<j$. 

We consider on each interval $[a,b]$ the uniform probability
distribution. Fix a parameter $d>0$. 
We generate the $n\times n$ matrix $R$ by induction on $n\geq 2$. 

\underline{$n=2$}: Pick $r_{1,2}$ randomly in $[0,d]$. 

\underline{$n\to n+1$}:
We define $r_{k,n+1}$ for $1\leq k\leq n$ by induction on $k$.

\underline{$k=1$}: Pick $r_{1,n+1}$ randomly in $[0,d]$.

\underline{$k\to k+1$}: Pick $r_{k+1,n+1}$ randomly in $[a,b]$ where
$$
  a:=\max_{1\leq i\leq k}|r_{i,k+1}-r_{i,n+1}|,\qquad 
  b:=\min_{1\leq i\leq k}|r_{i,k+1}+r_{i,n+1}|.  
$$

{\bf Products of metric spaces. }
For a real number $p\geq 1$ the $p$-product of two metric spaces is
described in terms of distance matrices as follows. For integers
$m,n\geq 2$ consider an $m\times m$ distance matrix $R$ and an
$n\times n$ distance matrix $S$. Their $p$-product $R\times_pS$ is the
$mn\times mn$ matrix $T$ with entries
$$
  t_{(i-1)n+k,(j-1)n+l} := (r_{i,j}^p + s_{k,l}^p)^{1/p},\qquad 1\leq
  i,j\leq m,\;1\leq k,l\leq n.  
$$
Note that this is again a distance matrix.

{\bf Computing the covering number. }
Let $R$ be an $n\times n$ distance matrix with parameter $d>0$ as
above. For a radius $r>0$, the covering number $\cov(R,r)$ is the
minimal number of closed $r$-balls covering the corresponding metric
space. Here is an algorithm to compute it. 

Define the $n\times n$ matrix $A$ with entries in $\{0,1\}$
$$
  a_{i,j} := \begin{cases} 1 & r_{i,j}\leq r, \cr 0 & r_{i,j}>r.
  \end{cases}
$$
Let $e$ be the $n$-vector with all entries $e_j=1$. Then
$$
  \cov(R,r) = \min\{e\cdot b \mid A\cdot b\geq e\}.
$$
Here the minimum is taken over all $n$-vectors with entries in
$\{0,1\}$, and $v\geq w$ for vectors means $v_j\geq w_j$ for all $j$. 

The interpretation is that we take the balls with centers $j$ for
$b_j=1$. The condition $A\cdot b\geq e$ means that these balls define
a covering, and $e\cdot b$ is the number of balls. 

Note that this computation of $\cov(R,r)$ is a problem in Integer
Linear Programming, for which there exist very efficient implementations.

\bibliographystyle{abbrv}
\bibliography{000_cov}

\end{document}